\documentclass{amsart}
\usepackage{amssymb, amsbsy, amsthm, amsmath, amstext, amsopn, verbatim}
\usepackage[all]{xy}
\usepackage{amsfonts}
\usepackage{amscd}
\usepackage{mathrsfs}
\usepackage{verbatim}
\usepackage{color}

\newtheorem{thm}{Theorem} [section]
\newtheorem{lemma}[thm]{Lemma}

\newtheorem{corollary}[thm]{Corollary}
\newtheorem{prop}[thm]{Proposition}
\newtheorem{notation}[thm]{Notation}

\newtheorem{assumptions}[thm]{Assumptions}
\newtheorem{assumption}[thm]{Assumption}

\theoremstyle{definition}

\newtheorem*{basic convention}{Basic Conventions}

\newtheorem{defn}[thm]{Definition}

\newtheorem*{HKK}{Hyperk\"ahler Kirwan Surjectivity Problem}

\theoremstyle{remark}

\newtheorem{remark}[thm]{Remark}

\begin{document}

\numberwithin{equation}{section}

\newcommand{\hs}{\mbox{\hspace{.4em}}}
\newcommand{\ds}{\displaystyle}
\newcommand{\bd}{\begin{displaymath}}
\newcommand{\ed}{\end{displaymath}}
\newcommand{\bcd}{\begin{CD}}
\newcommand{\ecd}{\end{CD}}

\newcommand{\proj}{\operatorname{Proj}}
\newcommand{\bproj}{\underline{\operatorname{Proj}}}
\newcommand{\spec}{\operatorname{Spec}}
\newcommand{\bspec}{\underline{\operatorname{Spec}}}
\newcommand{\pline}{{\mathbf P} ^1}
\newcommand{\pplane}{{\mathbf P}^2}
\newcommand{\coker}{{\operatorname{coker}}}
\newcommand{\ldb}{[[}
\newcommand{\rdb}{]]}

\newcommand{\Sym}{\operatorname{Sym}^{\bullet}}
\newcommand{\Symp}{\operatorname{Sym}}
\newcommand{\Pic}{\operatorname{Pic}}
\newcommand{\AAut}{\operatorname{Aut}}
\newcommand{\PAut}{\operatorname{PAut}}

\newcommand{\too}{\twoheadrightarrow}
\newcommand{\C}{{\mathbf C}}
\newcommand{\cA}{{\mathcal A}}
\newcommand{\cS}{{\mathcal S}}
\newcommand{\cV}{{\mathcal V}}
\newcommand{\cM}{{\mathcal M}}
\newcommand{\bA}{{\mathbf A}}
\newcommand{\aline}{\mathbb{A}^1}
\newcommand{\cB}{{\mathcal B}}
\newcommand{\cC}{{\mathcal C}}
\newcommand{\cD}{{\mathcal D}}
\newcommand{\D}{{\mathcal D}}
\newcommand{\cs}{{\mathbf C} ^*}
\newcommand{\boldc}{{\mathbf C}}
\newcommand{\cE}{{\mathcal E}}
\newcommand{\cF}{{\mathcal F}}
\newcommand{\cG}{{\mathcal G}}
\newcommand{\G}{{\mathbf G}}
\newcommand{\fg}{{\mathfrak g}}
\newcommand{\ft}{\mathfrak t}
\newcommand{\bH}{{\mathbf H}}
\newcommand{\cH}{{\mathcal H}}
\newcommand{\cI}{{\mathcal I}}
\newcommand{\cJ}{{\mathcal J}}
\newcommand{\cK}{{\mathcal K}}
\newcommand{\cL}{{\mathcal L}}
\newcommand{\baL}{{\overline{\mathcal L}}}
\newcommand{\M}{{\mathcal M}}
\newcommand{\bM}{{\mathbf M}}
\newcommand{\bm}{{\mathbf m}}
\newcommand{\cN}{{\mathcal N}}
\newcommand{\theo}{\mathcal{O}}
\newcommand{\cP}{{\mathcal P}}
\newcommand{\cR}{{\mathcal R}}
\newcommand{\boldp}{{\mathbf P}}
\newcommand{\boldq}{{\mathbf Q}}
\newcommand{\bbL}{{\mathbf L}}
\newcommand{\cQ}{{\mathcal Q}}
\newcommand{\cO}{{\mathcal O}}
\newcommand{\Oo}{{\mathcal O}}
\newcommand{\OX}{{\Oo_X}}
\newcommand{\OY}{{\Oo_Y}}
\newcommand{\otY}{{\underset{\OY}{\ot}}}
\newcommand{\otX}{{\underset{\OX}{\ot}}}
\newcommand{\cU}{{\mathcal U}}
\newcommand{\cX}{{\mathcal X}}
\newcommand{\cW}{{\mathcal W}}
\newcommand{\boldz}{{\mathbf Z}}
\newcommand{\cZ}{{\mathcal Z}}
\newcommand{\qgr}{\operatorname{qgr}}
\newcommand{\gr}{\operatorname{gr}}
\newcommand{\coh}{\operatorname{coh}}
\newcommand{\End}{\operatorname{End}}
\newcommand{\Hom}{\operatorname{Hom}}
\newcommand{\uHom}{\underline{\operatorname{Hom}}}
\newcommand{\uHomY}{\uHom_{\OY}}
\newcommand{\uHomX}{\uHom_{\OX}}
\newcommand{\Ext}{\operatorname{Ext}}
\newcommand{\bExt}{\operatorname{\bf{Ext}}}
\newcommand{\Tor}{\operatorname{Tor}}

\newcommand{\inv}{^{-1}}
\newcommand{\airtilde}{\widetilde{\hspace{.5em}}}
\newcommand{\airhat}{\widehat{\hspace{.5em}}}
\newcommand{\nt}{^{\circ}}
\newcommand{\del}{\partial}

\newcommand{\supp}{\operatorname{supp}}
\newcommand{\GK}{\operatorname{GK-dim}}
\newcommand{\hd}{\operatorname{hd}}
\newcommand{\id}{\operatorname{id}}
\newcommand{\res}{\operatorname{res}}
\newcommand{\lrar}{\leadsto}
\newcommand{\im}{\operatorname{Im}}
\newcommand{\HH}{\operatorname{H}}
\newcommand{\TF}{\operatorname{TF}}
\newcommand{\Bun}{\operatorname{Bun}}
\newcommand{\Hilb}{\operatorname{Hilb}}
\newcommand{\Fact}{\operatorname{Fact}}
\newcommand{\F}{\mathcal{F}}
\newcommand{\nthord}{^{(n)}}
\newcommand{\Aut}{\underline{\operatorname{Aut}}}
\newcommand{\Gr}{\operatorname{Gr}}
\newcommand{\Fr}{\operatorname{Fr}}
\newcommand{\GL}{\operatorname{GL}}
\newcommand{\gl}{\mathfrak{gl}}
\newcommand{\SL}{\operatorname{SL}}
\newcommand{\ff}{\footnote}
\newcommand{\ot}{\otimes}
\def\Ext{\operatorname {Ext}}
\def\Hom{\operatorname {Hom}}
\def\Ind{\operatorname {Ind}}
\def\bbZ{{\mathbb Z}}

\newcommand{\nc}{\newcommand}
\newcommand{\on}{\operatorname}
\nc{\cont}{\on{cont}}
\nc{\rmod}{\on{mod}}
\nc{\Mtil}{\widetilde{M}}
\nc{\wb}{\overline}
\nc{\wt}{\widetilde}
\nc{\wh}{\widehat}
\nc{\sm}{\setminus}
\nc{\mc}{\mathcal}
\nc{\mbb}{\mathbb}
\nc{\Mbar}{\wb{M}}
\nc{\Nbar}{\wb{N}}
\nc{\Mhat}{\wh{M}}
\nc{\pihat}{\wh{\pi}}
\nc{\JYX}{\cJ_{Y\leftarrow X}}
\nc{\phitil}{\wt{\phi}}
\nc{\Qbar}{\wb{Q}}
\nc{\DYX}{\D_{Y\leftarrow X}}
\nc{\DXY}{\D_{X\to Y}}
\nc{\dR}{\stackrel{\bbL}{\underset{\D_X}{\ot}}}
\nc{\Winfi}{\cW_{1+\infty}}
\nc{\K}{{\mc K}}
\nc{\unit}{{\bf \on{unit}}}
\nc{\boxt}{\boxtimes}
\nc{\xarr}{\stackrel{\rightarrow}{x}}
\nc{\Cnatbar}{\overline{C}^{\natural}}
\nc{\oJac}{\overline{\on{Jac}}}
\nc{\gm}{{\mathbf G}_m}
\nc{\Loc}{\on{Loc}}
\nc{\Bm}{\operatorname{Bimod}}
\nc{\lie}{{\mathfrak g}}
\nc{\lb}{{\mathfrak b}}
\nc{\lien}{{\mathfrak n}}
\nc{\e}{\epsilon}
\nc{\eu}{\mathsf{eu}}
\nc{\bs}{\backslash}

\nc{\Gm}{{\mathbb G}_m}
\nc{\Gabar}{\wb{\G}_a}
\nc{\Gmbar}{\wb{\G}_m}
\nc{\PD}{{\mathbb P}_{\D}}
\nc{\Pbul}{P_{\bullet}}
\nc{\PDl}{{\mathbb P}_{\D(\lambda)}}
\nc{\PLoc}{\mathsf{MLoc}}
\nc{\Tors}{\on{Tors}}
\nc{\PS}{{\mathsf{PS}}}
\nc{\PB}{{\mathsf{MB}}}
\nc{\Pb}{{\underline{\operatorname{MBun}}}}
\nc{\Ht}{\mathsf{H}}
\nc{\bbH}{\mathbb H}
\nc{\gen}{^\circ}
\nc{\Jac}{\operatorname{Jac}}
\nc{\sP}{\mathsf{P}}
\nc{\sT}{\mathsf{T}}
\nc{\bP}{{\mathbb P}}
\nc{\otc}{^{\otimes c}}
\nc{\Det}{\mathsf{det}}
\nc{\PL}{\on{ML}}
\nc{\sL}{\mathsf{L}}

\nc{\ml}{{\mathcal S}}
\nc{\Xc}{X_{\on{con}}}
\nc{\Z}{{\mathbb Z}}
\nc{\resol}{\mathfrak{X}}
\nc{\map}{\mathsf{f}}
\nc{\gK}{\mathbb{K}}
\nc{\bigvar}{\mathsf{W}}
\nc{\Dhol}{\mathcal{D}_{\on{hol}}}
\nc{\Dcrh}{\mathcal{D}_{\on{crh}}}
\nc{\fX}{\mathfrak{X}}
\nc{\A}{\mathsf{A}}
\nc{\vD}{\mathbb{D}}
\nc{\openset}{U}
\nc{\cat}{\mathcal{C}}
\nc{\jb}{j_{\beta}}
\nc{\cY}{\mathcal{Y}}

\title[]{Morse Decomposition for $\D$-Module Categories on Stacks}
\author{Kevin McGerty}
\address{Mathematical Institute\\University of Oxford\\Oxford OX1 3LB, UK}
\email{mcgerty@maths.ox.ac.uk}
\author{Thomas Nevins}
\address{Department of Mathematics\\University of Illinois at Urbana-Champaign\\Urbana, IL 61801 USA}
\email{nevins@illinois.edu}

\begin{abstract}
Let $\mathcal{Y}$ be a smooth algebraic stack exhausted by quotient stacks.  Given a Kirwan-Ness stratification of the cotangent stack $T^*\mathcal{Y}$, we establish a recollement package for twisted $\D$-modules on $\mathcal{Y}$, gluing the category from subquotients described via modules microsupported on the Kirwan-Ness strata of $T^*\mathcal{Y}$.   The package includes unusual existence and ``preservation-of-finiteness'' properties for  functors of the {\em full} category of twisted $\D$-modules, extending the standard functorialities for holonomic modules.   In the case that $\mathcal{Y} = X/G$ is a quotient stack, our results provide a higher categorical analogue of the Atiyah-Bott--Kirwan--Ness ``equivariant perfection of Morse theory'' for the norm-squared of a real moment map.  As a consequence, we deduce a modified form of Kirwan surjectivity for the cohomology of hyperk\"ahler/algebraic symplectic quotients of cotangent bundles.  
\end{abstract}

\maketitle

\section{Introduction}
Let $X$ be a smooth complex algebraic variety with the action of a reductive group $G$.  Choosing a semistable locus $X^{ss}\subset X$ in the sense of geometric invariant theory (GIT),  work of many authors  (cf. \cite{AB, Kirwan, Ness}) has closely linked the cohomology of $X^{ss}/G$ to the equivariant cohomology of $X$ via a Morse stratification, commonly called the Kirwan-Ness stratification, of $X^u = X\smallsetminus X^{ss}$.  Recent work (cf. \cite{JKK, DWWW}) addresses the cohomology of hyperk\"ahler quotients via Morse theory.  

\vspace{.5em}

In the present paper, we establish a strong form of Morse decomposition one categorical level higher, for (twisted) equivariant $\D$-modules on $X$, or more generally $\D$-modules on stacks.  We deduce cohomological consequences for quotients in hyperk\"ahler and algebraic symplectic geometry.

\subsection{$\D$-Modules and Microlocal KN Stratifications}
The main object of study in this paper is the unbounded derived category (or rather its dg enhancement) $\D(X/G)$ of quasicoherent $\D$-modules on the quotient stack $X/G$, or, equivalently, an appropriately defined derived category of $G$-equivariant $\D$-modules on $X$; see Definition \ref{def:D-mod} for details.\ff{We note that this category is often not the derived category of an abelian category.  Accordingly, our notation for functors between such categories omits any of the modifiers such as $\mathbb{L}, \mathbb{R}$ typically used for derived functors.}
  More generally, associated to a quantum comoment map
$\mu_c: \mathfrak{g}\rightarrow \D(X)$, one can consider a twisted form $\D(X/G,c)$ of $\D(X/G)$, as explained in \cite{McNderived} or \cite{McN} in notation consistent with what we use here.  More generally still, we let $\mathcal{Y}$ denote an algebraic stack that has a Zariski-open cover by stacks of the form $X/G$ where $X$ is a smooth complex variety and $G$ is a reductive algebraic group; we say such a stack is {\em exhausted by (smooth) quotient stacks}; we write $\D(\mathcal{Y},c)$ for the derived category of $c$-twisted quasi-coherent $\D$-modules on $\mathcal{Y}$.   In general, this category is defined via a homotopy limit; in the special case when $\mathcal{Y}=X/G$ and the classical moment map $T^*X\rightarrow\mathfrak{g}^*$ is flat, the category is the (dg enhanced) derived category of the abelian category of twisted $G$-equivariant $\D$-modules on $X$.

Letting $X$ be a smooth variety with $G$-action, consider any  $G$-stable conical open subset $U\subseteq T^*X$ and the complementary closed subset $K\subset T^*X$.  Let $\D(X/G,c)_K$ denote the subcategory of $\D(X,G,c)$ consisting of objects whose pullback to $X$ has cohomologies microsupported in $K$.  The {\em microlocal category} associated to $U$ is then defined to be $\D(U/\!\!/G,c):=\D(X/G,c)/\D(X/G,c)_K$.  This category can sometimes be understood as a derived category of modules for an algebra of twisted differential operators microrestricted to the open subset $U/\!\!/G$ of the cotangent stack $T^*(X/G)$ (see Section 4.5 of \cite{McNderived} for a discussion).  A similar definition applies for a conical closed subset $\mathcal{K}\subseteq T^*\mathcal{Y}$ for a stack $\mathcal{Y}$ with open complement $\cU\subseteq T^*\mathcal{Y}$ to give categories
$\D(\mathcal{Y},c)_{\mathcal{K}}$ and $\D(\cU,c)$.  

We assume that $T^*\mathcal{Y}$ comes equipped with a {\em Kirwan-Ness (KN) stratification}: this is a decomposition into smooth, locally closed substacks with nice properties, which, in the case $\mathcal{Y} = X/G$, reflect the equivariant Morse theory (for an appropriate Morse function) of $T^*X$.  See  Section \ref{sec:KN} for details of KN stratifications and \cite{Kirwan} for the connection to Morse theory.  Our definition of a KN stratification for $T^*\mathcal{Y}$ requires that each stratum $\mathcal{K}$ arises, on a chart $X/G\subseteq \mathcal{Y}$, as $(\mu\inv(0)\cap K)/G$ where $K\subseteq T^*X$ is part of a KN stratification and $\mu: T^*X\rightarrow \mathfrak{g}^*$ is the classical moment map.  

\subsection{Categorical Morse Decomposition}\label{sec:main-results}
Let $\mathcal{K}\subseteq \mathcal{K}'\subset T^*\mathcal{Y}$ be any two closed unions of KN strata, and $\cU = T^*\mathcal{Y}\smallsetminus \mathcal{K}$ and $\cU^\circ = T^*\mathcal{Y}\smallsetminus \mathcal{K}'$ the complementary open sets, so $j: \cU^\circ\hookrightarrow \cU$ is an open immersion.  Let $i: Z:= \mathcal{K}'\cap \cU\hookrightarrow \cU$ denote the closed immersion, an inclusion of a union of KN strata.    There is a natural $t$-exact and continuous functor of restriction,
\begin{equation}\label{restrictions}
j^*: \D(\cU,c)\rightarrow \D(\cU^\circ,c).
\end{equation}
This functor admits a right adjoint $j_*: \D(\cU^\circ,c)\rightarrow \D(\cU,c)$.  There is a subcategory $i_*: \D(\cU,c)_Z \hookrightarrow \D(\cU,c)$ of objects microsupported on $Z$.  The functor $i_*$ has a right adjoint $i^!$.  

Each of the categories $\D(\cU,c)$, etc. has two important full subcategories: the category $\D(\cU,c)^c$ of compact objects and the category $\D(\cU,c)_{\on{coh}}$ of bounded complexes with coherent (as $\D$-modules) cohomologies (or just {\em coherent objects}).  We have $\D(\cU,c)^c\subseteq \D(\cU,c)_{\on{coh}}$, but the characterization of $\D(\cU,c)^c$ is, in general, subtle.

The main results of the paper are the following.
\begin{thm}\label{main theorem}
\mbox{}
\begin{enumerate}
\item The functor $j^*$ admits a left adjoint $j_!$.  
\item The functor $i_*$ admits a left adjoint $i^*$.  
\item The functors in the adjoint triples $(j_!, j^*, j_*)$, $(i^*, i_*, i^!)$ preserve  compactness and coherence.
%\item The left adjoint $j_!$ takes compact objects to compact objects.
%\item The right adjoint $j_*$ of $j^*$ takes bounded coherent complexes to bounded coherent complexes.
%\item The functor $i^!$ admits a continuous right adjoint.
%\item The functor $i_*$ admits a left adjoint $i^*$.
%\item The left adjoint $i^*$ takes compact objects to compact objects.  
%\item The right adjoint $i^!$ takes bounded coherent complexes to bounded coherent complexes.

\end{enumerate}
\end{thm}
\begin{thm}\label{second main thm}
Let $Z = \cU\smallsetminus \cU^\circ$.  The adjoint triples of functors $(i^*,i_*,i^!)$ and $(j_!,j^*,j_*)$ assemble into a {\em recollement diagram},  
\bd
\xymatrix{
\D(\cU,c)_Z  \ar[r]^{i_*}  & 
\D(\cU,c) \ar@/^1pc/[l]^{i^!}\ar@/_1pc/[l]_{i^*} \ar[r]^{j^*}  
& \D(\cU^\circ,c) \ar@/^1pc/[l]^{j_*}\ar@/_1pc/[l]_{j_!}.
}
\ed

That is, there are exact triangles of functors
\bd
i_*i^! \rightarrow \on{Id} \rightarrow j_*j^* \xrightarrow{[1]}\hspace{2em} \text{and}\hspace{2em}  j_!j^*\rightarrow \on{Id}\rightarrow i_*i^*\xrightarrow{[1]},
\ed
so that $j^*\circ i_*\simeq 0$, and the functors $j_!, j_*$, and $i_*$ are full embeddings.  
\end{thm}

\vspace{.5em}

\noindent
By Theorem \ref{main theorem} and formal properties, analogous recollements hold for $\D(\cU,c)^c$ and $\D(\cU,c)_{\on{coh}}$.

\vspace{.5em}

The results immediately reduce, in the case $\mathcal{Y}=X/G$, to more classical statements about equivariant derived categories.  We write $K\subseteq K'\subseteq T^*X$ for the closed unions of KN strata in $T^*X$ and  $U = T^*X\smallsetminus K$, $U^\circ = T^*X\smallsetminus K'$ for the complements.
\begin{corollary}
Let $Z = U\smallsetminus U^\circ$.  The adjoint triples of functors $(i^*,i_*,i^!)$ and $(j_!,j^*,j_*)$ assemble into a recollement diagram,  
\bd
\xymatrix{
\D(U/\!\!/G,c)_Z  \ar[r]^{i_*}  & 
\D(U/\!\!/G,c) \ar@/^1pc/[l]^{i^!}\ar@/_1pc/[l]_{i^*} \ar[r]^{j^*}  
& \D(U^\circ/\!\!/G,c) \ar@/^1pc/[l]^{j_*}\ar@/_1pc/[l]_{j_!}.
}
\ed
That is, there are adjoint triples $(i^*,i_*,i^!)$, $(j_!,j^*,j_*)$ and exact triangles of functors
\bd
i_*i^! \rightarrow \on{Id} \rightarrow j_*j^* \xrightarrow{[1]}\hspace{2em} \text{and}\hspace{2em}  j_!j^*\rightarrow \on{Id}\rightarrow i_*i^*\xrightarrow{[1]},
\ed
so that $j^*\circ i_*\simeq 0$, and the functors $j_!, j_*$, and $i_*$ are full embeddings.  These functors satisfy all the properties asserted in Theorem \ref{main theorem}.
\end{corollary}
An analogous assertion to the existence of $j_!$ in Theorem \ref{main theorem} (in the special case when $K$ is empty and $K'$ is the entire unstable locus)
for deformation quantization modules appears as \cite[Lemma~5.18]{BPW}, though the method of attack is completely different than the one here.   

The existence of  a {\em left} adjoint $j_!$ to $j^*$ is a somewhat surprising feature of the equivariant quantum world, i.e., of (twisted) $\D$-modules.\ff{We expect analogues to hold for equivariant sheaves of modules over deformation quantizations of smooth algebraic symplectic varieties, via essentially similar techniques.}  We note however that a similar collection of statements {\em never} holds (except in trivial cases) if we replace the category of $\D$-modules (or DQ modules) with quasicoherent sheaves.  It also essentially never exists for the {\em full} category of $\D$-modules if $U\hookrightarrow X$ is a nonempty, open, proper subset of an ordinary variety.  On the other hand, the functor $j_!$ is always defined for {\em holonomic} $\D$-modules---but only for reasons special to holonomic modules.  What is different in our setting is that the particular equivariant (or Morse-theoretic) structure of the KN strata has very strong consequences for pushforwards across strata.  Similarly, the preservation-of-coherence statement of part (3) of the theorem never holds in the quasicoherent (commutative) world; 
however, there is an analogue to $j_*(\D(U/\!\!/G)_{\on{coh}})\subset \D(X/G)_{\on{coh}}$, namely the semiorthogonal complement to the kernel of $j^*$  in the beautiful work of Halpern-Leistner \cite{HL}.\ff{We emphasize that in the world of $\D$-modules this semiorthogonal decomposition is canonically determined by adjunction once one fixes a KN stratification, whereas in the coherent world it depends on a further choice of grade restriction windows.}

Our main results incarnate, one categorical level higher, the equivariantly perfect Morse decompositions of \cite{AB,Kirwan, Ness} for the Morse flow associated to the norm-squared of the real moment map, associated to the action of a compact real form of $G$ on $T^*X$, under the choice of a K\"ahler metric.  
The comprehensive reference \cite{BPW} convincingly advocates the view that $j^*$ should be considered as a ``categorical Kirwan map.''  
The general paradigm of  categorical Morse decomposition is eloquently introduced in
\cite{Nadler}, where results are established in real symplectic geometry (of the Fukaya category).  We note that \cite{Nadler} considers the ``classical'' Fukaya category of Lagrangian branes, and thus the relevant functoriality is closer kin to that of holonomic $\D$-modules than to the full category of $\D$-modules (analogues of general coisotropic branes) that we consider here.  A ``categorical Morse theory'' in the spirit of Bia\l ynicki-Birula (but for non-Hamiltonian $\Gm-$actions) is explored in \cite{BDMN}.     

\subsection{Kirwan Surjectivity}
The existence of adjoints in Theorem \ref{main theorem} has strong consequences for cohomology of quotients.  Namely, suppose now that $X$ is a smooth, connected, quasi-compact scheme with action by a reductive $G$.  Assume in addition that the moment map $\mu: T^*X\rightarrow \mathfrak{g}^*$ is flat, so that $\mu\inv(0)$ is a complete intersection in $T^*X$.  Assume that $T^*X$ is equipped with a Kirwan-Ness stratification, and suppose that the $G$-action on the intersection $\mu\inv(0)^{ss}:=\mu\inv(0)\cap T^*X^{ss}$ of the semistable (open) stratum with $\mu\inv(0)$ consists of free $G$-orbits, so that there is a smooth scheme $\mu\inv(0)^{ss}/G$ for which the projection $\mu\inv(0)^{ss}\rightarrow \mu\inv(0)^{ss}/G$ is a principal $G$-bundle.  

Recall \cite{JKK} that if the cotangent bundle $T^*X$ admits a hyperk\"ahler metric invariant under the compact real form $G_c$ of $G$, then $\mu\inv(0)^{ss}/G$ can typically be interpreted as the quotient by $G_c$ of a fiber of a hyperk\"ahler moment map.  The natural map
\begin{equation}\label{eq:HKKirwan}
H^*_{dR}(X/G) \cong H^*_G(T^*X)\longrightarrow H^*_G\big(\mu\inv(0)^{ss}\big) = H^*_{dR}(\mu\inv(0)^{ss}/G)
\end{equation}
is the {\em hyperk\"ahler Kirwan map} \cite{JKK}.   
\begin{HKK}
When is the map \eqref{eq:HKKirwan} surjective?
\end{HKK}
\noindent
Let $I(X)$ denote the {\em inertia scheme} for the $G$-action on $X$: that is the fiber product of the diagram $G\times X\xrightarrow{(a,p_2)} X\times X \xleftarrow{\Delta} X$, where $(a,p_2)$ denotes the product of the action map $a:G\times X\rightarrow X$ and the projection on the second factor.  The quotient $I(X)/G$ is the inertia stack of $X/G$.  Then:
\begin{thm}\label{thm:Kirwan}
There is a natural surjective homomorphism
\begin{equation}\label{eq:Kirwan}
H^*_{dR}\big(I(X)/G\big) = H^*_{G}\big(I(X)\big)\longrightarrow H^*_{dR}(\mu\inv(0)^{ss}/G).
\end{equation}
\end{thm}
\noindent
We note that Theorem \ref{thm:Kirwan} does {\em not} require a hyperk\"ahler metric, only the structure of a KN stratification of $T^*X$ in algebraic symplectic geometry.  The theorem is best interpreted as a modified form of hyperk\"ahler Kirwan surjectivity, since Theorem \ref{thm:Kirwan} proves that \eqref{eq:Kirwan} is {\em always surjective}, which is not true for the Kirwan map itself.  The theorem may be interpreted as saying that a Kirwan map incorporating {\em twisted sectors}, i.e., the cohomology of the full inertia stack, is a more tractable object of study.  

N. Proudfoot and B. Webster have explained to us that one can use Theorem \ref{thm:Kirwan} and cohomological purity arguments to prove the surjectivity of the map \eqref{eq:HKKirwan} for $X$ a linear representation of a connected reductive group $G$ (in particular, for quiver varieties).  Details will appear in an appendix authored by Proudfoot and Webster to a future version of this paper.

A particularly interesting example, not specifically addressed here, is the (non--quasi-compact) moduli stack $\on{Bun}_G(C)$ of $G$-bundles for a smooth projective curve $C$. Here the inertia stack is the moduli stack of pairs consisting of a $G$-bundle and an automorphism of it, sometimes called a {\em group-like Hitchin space}; the target of the map corresponding to \eqref{eq:Kirwan} is the cohomology of the moduli space of (semi)stable Higgs bundles.
We hope to carry out a more detailed analysis of this example, including a comparison to \cite{DWWW}, elsewhere.

Theorem \ref{thm:Kirwan} is, modulo known results, a formal consequence of Theorem \ref{main theorem} as interpreted in Hochschild homology, and is explained in
Section \ref{sec:Kirwan}.  

\subsection{Antecedents and Further Directions}
Parts of Theorem \ref{main theorem}  have antecedents in the literature.  We have already mentioned \cite{BPW}.  However, a more direct technical inspiration for the present paper, and the source of several of the categorical ideas used in the proof, is the work of Drinfeld--Gaitsgory proving compact generation of the category of $\D$-modules on $\on{Bun}_G(C)$.  The paper \cite{DG} develops some of the assertions of Theorem \ref{main theorem} for a particular stratification of the stack $\on{Bun}_G(C)$.  
Theorem \ref{main theorem} springs from the realization that the statement of, and geometry behind, the vanishing theorem of \cite{McN} could be combined with suitably generalized versions of the more categorical parts of \cite{DG} to prove analogous statements in the generality of the present paper.  We note that Theorem \ref{second main thm} represents a particularly illuminating re-packaging of the structure already present when one knows the assertions of Theorem \ref{main theorem}.

In \cite{DG}, the existence of $j_!$ and its preservation of compactness are used to produce compact objects;
in particular, Theorem \ref{main theorem} together with the basic theory of $\D$-modules on non--quasi-compact stacks from \cite{DG} provides a new proof of the following:
\begin{thm}[Drinfeld-Gaitsgory, \cite{DG}]\label{cor:cpt gen}
Let $C$ be a smooth projective curve, $G$ a reductive group,  $c\in\C$.  The unbounded derived category of  quasi-coherent $\D(\mathsf{det}^{\otimes c})$-modules on $\on{Bun}_G(C)$ is compactly generated.
\end{thm}
The present paper, however, does {\em not} include full details of a proof of Theorem \ref{cor:cpt gen}: rather, we intend to pursue compact generation questions, including the more subtle problem of compact generation of the ind-holonomic category (cf. \cite{BNa}), elsewhere.   The solutions to the latter problem hinges on a more precise structure theory for the subquotients of modules supported on KN strata, which we plan to develop in forthcoming work.

\subsection*{Acknowledgments}
We are grateful to David Ben-Zvi for extensive consultation and comments on an ealier draft, and to David Nadler, Charles Rezk, and Nick Rozenblyum for helpful conversations. 

This work was conceived while the second author was supported by a Visiting Fellowship at All Souls College, Oxford; we are grateful to the college for providing this support for our collaboration.  Parts of this work were also carried out at the Aspen Center for Physics; we are grateful to the ACP for its hospitality and to the NSF, which, via grant PHYS-1066293, supports the ACP.  
The first author was supported by a 
 Royal Society research fellowship.  The second author was supported by NSF grant DMS-1159468 and NSA grant H98230-12-1-0216.  

\subsection*{Conventions}
We work throughout with pre-triangulated dg categories over a fixed algebraically closed field of characteristic zero (which we abusively denote $\C$).  Basic properties can be found in \cite{Drinfeld, Keller2, To, Gaitsgory}.  Given a dg category $C$ with $t$-structure, 
we write $C^{\heartsuit}$ for the abelian category that is the heart of the $t$-structure.  Given a functor of dg categories, ``faithful'' will always mean ``quasi-faithful,'' ``equivalence'' means ``quasi-equivalence,'' etc.

\section{KN Stratifications and Equivariant Geometry Near a Stratum}\label{sec:KN}
This section lays out basics of KN stratifications and reviews features of the symplectic geometry near a KN stratum, following \cite{McN}, that play a crucial role in this paper.
\subsection{Kirwan-Ness Stratifications}\label{sec:KN-for-vars}
We begin by reviewing KN stratifications and discussing a reasonable extension to cotangent stacks of algebraic stacks exhausted by global quotients.
\subsubsection{KN Stratification of Varieties}
Let $G$ be a reductive group and $\sT$ a maximal torus of $G$, with $W$ the corresponding Weyl group. Let $Y(\sT) = \on{Hom}(\mathbb G_m,\sT)$ be the group of $1$-parameter subgroups of $\sT$, and $X(\sT) = \on{Hom}(\sT,\mathbb G_m)$, the group of characters. Let $Y_\mathbb Q = Y(\sT)_\mathbb Q$ be the $\mathbb Q$-vector space\footnote{Since $\sT$ is a torus, $Y_\mathbb Q$ is a rational form of $\on{Lie}(\sT)$. For the reductive group $G$, it makes sense to define the set of ``rational $1$-parameter subgroups'' but it no longer has the structure of a $\mathbb Q$-vector space.}  $Y(\sT)\otimes_\mathbb Z \mathbb Q$ and $X_\mathbb Q = X(\sT)\otimes_\mathbb Z \mathbb Q$ similarly. Let $q\colon Y_{\mathbb Q} \to \mathbb Q$ denote a $W$-invariant, integral, positive definite quadratic form, and $d$ the induced metric. The quadratic form allows us to identify $X_\mathbb Q$ and $Y_{\mathbb Q}$, which we do henceforth. 

Suppose $\bigvar$ is a smooth $G$-variety equipped with a $G$-equivariant line bundle $\mathscr{L}$.  Let $\mathsf{KN} = \{\beta\}$ be a finite collection of 1-parameter subgroups of $\sT$.  
We define a partial order $<$ on $\mathsf{KN}$ by setting $\beta < \beta'$ if $q(\beta)<q(\beta')$. 
Given $\beta\in\mathsf{KN}$, write $Z_\beta = \bigvar^{\beta(\Gm)}$ for the fixed-point locus of $\beta$.  Let 
\bd
\ds Y_\beta = \{x\in \bigvar \;|\; \lim_{t\rightarrow 0} \beta(t)\cdot x \in Z_\beta\}, \;\; \text{and write}\;\;
\on{pr}_\beta: Y_\beta\rightarrow Z_\beta
\ed
for the corresponding projection.  Write $\bigvar^{ss}$ for the open subset of $\bigvar$ complementary to the common zero locus of $G$-invariant elements of $H^0(\bigvar, \mathscr{L}^N)$ for all $N>0$.  

Let $L_\beta \subseteq G$ denote the centralizer of $\beta$ in $G$, and let $P_\beta$ denote the parabolic subgroup of $G$ whose Lie algebra is spanned by the nonnegative $\beta$-weight spaces in $\mathfrak{g}$; in particular, 
$L_\beta$ is the Levi factor of $P_\beta$.  Then $L_\beta$ preserves $Z_\beta$, $P_\beta$ preserves $Y_\beta$, and $\on{pr}_\beta$ is $P_\beta$-equivariant where $P_\beta$ acts on $Z_\beta$ via $P_\beta\rightarrow L_\beta$.  
Moreover, for each connected component $Z_{\beta, i}$ of $Z_\beta$ and $x\in Z_{\beta, i}$, $L_\beta$ acts via a character $\lambda_{\beta, i}: L_\beta \rightarrow \Gm$ on the fiber $\mathscr{L}(x)$ of $\mathscr{L}$ over $x$ (and this character does not depend on the choice of $x\in Z_{\beta, i}$).

\begin{defn}\label{def:KN-strat}
A {\em KN stratification} of $\bigvar$ with respect to $\mathsf{KN}$ consists of a choice, 
 for each $\beta$, of an $L_\beta$-stable open subset $Z_\beta^{ss}\subseteq Z_\beta$ satisfying:
 \begin{enumerate} 
\item for each $\beta$ and each component $Z_{\beta, i}$ , the complement $Z_{\beta, i}\smallsetminus Z_\beta^{ss}$ in $Z_{\beta, i}$ is cut out by a collection of $\lambda_{\beta, i}$-semi-invariant elements of $H^0\big(Z_{\beta, i}, \mathscr{L}^N|_{Z_{\beta, i}}\big)$ for some $N \gg 0$.  
\item Defining $Y_\beta^{ss} = \on{pr}_\beta\inv(Z_\beta^{ss})$, $Y_\beta^{ss}$ is a $\beta$-equivariant affine bundle over $Z_\beta^{ss}$ via $\on{pr}_\beta$.
\item Letting $S_\beta = G\cdot Y_\beta^{ss}$, we have $S_\beta \cong G\times_{P_\beta} Y_\beta^{ss}$.
\item  The collection of subsets $\{\bigvar^{ss}\}\cup \{S_\beta\; |\; \beta\in\mathsf{KN}\}$ stratifies $\bigvar$:
\begin{enumerate}
\item The stratum closures $\overline{S}_{\beta}$ satisfy $\overline{S}_{\beta}\subseteq \bigcup_{\beta'\geq \beta} S_{\beta'}$.\ff{It is \textit{not} assumed to be the case that the closure of a stratum is a union of strata.}
\item  $\bigvar = \bigvar^{ss}\coprod \big(\coprod_\beta S_\beta\big)$.
\end{enumerate}
\end{enumerate}
\end{defn}

\subsubsection{KN Stratifications of Cotangent Stacks} 
The results of this paper extend immediately to categories of twisted $\D$-modules on stacks that are exhausted by global quotient stacks $X/G$ where $X$ is a smooth variety and $G$ is a reductive algebraic group.  A typical example is the stack $\on{Bun}_H(C)$ of principal $H$-bundles on a curve $C$, where $H$ is a reductive group.

Suppose $\mathcal{Y}$ is such a stack with covering by $X_i/G_i$.  For each $i$, the cotangent stack $T^*(X_i/G_i)$, which is a Zariski-open substack of $T^*\mathcal{Y}$, is canonically isomorphic to the quotient $\mu_i\inv(0)/G_i$, where $\mu_i: T^*X_i \rightarrow \mathfrak{g}_i$ is the canonical classical moment map.  Given a KN stratification $\{S_\beta\}$ of $T^*X_i$ as in Section \ref{sec:KN-for-vars}, we can form a decomposition of $T^*(X_i/G_i)$ into locally closed substacks $\mathcal{S}_\beta$ by defining 
$\mathcal{S}_\beta = (\mu_i\inv(0)\cap S_\beta)/G_i$.  We call this the {\em induced} decomposition of $T^*(X_i/G_i)$.

\begin{defn}
A {\em KN stratification} $T^*\mathcal{Y} = \coprod \mathcal{S}_\beta$ for a stack $\mathcal{Y}$ exhausted by smooth quotient stacks $X_i/G_i$ is a decomposition into locally closed substacks such that, for each $i$, the intersections $\mathcal{S}_\beta\cap T^*(X_i/G_i)$ can be identified with the induced decomposition of $T^*(X_i/G_i)$ associated to a KN stratification of $T^*X_i$ as above.
\end{defn}

\subsection{Equivariant Geometry Near a KN Stratum}\label{sec:sympl-geom}
This section largely summarizes some symplectic geometry from \cite{McN}.

A KN stratum $S_\beta$ is labelled by (the Weyl group orbit of) a 1-parameter subgroup $\beta: \Gm\rightarrow \sT\subseteq G$ in a fixed maximal torus $\sT$ of $G$.  As in Section 12 of \cite{Kirwan}, $\beta$ determines a parabolic subgroup $P_\beta$ of $G$: letting $L_\beta$ denote the centralizer of $\beta$ in $G$, we let $P_\beta$ denote the subgroup whose Lie algebra is spanned by $\on{Lie}(L_\beta)$ and the positive $\beta$-weight subspaces in $\mathfrak{g}$.  The sum of positive $\beta$-weight subspaces is a nilpotent Lie sub-algebra $\mathfrak{n}\subset \mathfrak{g}$; we let $\mathfrak{n}^-$ denote the opposite nilpotent subalgebra, the sum of negative $\beta$-weight subspaces in $\mathfrak{g}$, and let $U^- = U_{P_\beta}^-\subset G$ denote the corresponding unipotent subgroup of $G$. 

We will write 
$\gK = U^-\ltimes \Gm$,
where $\Gm$ acts on $U^-$ via $\beta$ and the adjoint action of $G$.
The group $\gK$ acts naturally on $T^*X$ via $G$. We view $U^-$ as a $\gK$-variety where $\beta(\mathbb G_m)$ acts by conjugation and $U^-$ by left translation. 

\begin{lemma}
\label{action map}
The action map $a\colon U^-\times Y_\beta \rightarrow T^*W$ is a $\mathbb K$-equivariant bijection onto an open dense subset of $S_\beta = G\cdot Y_\beta \cong G\times_{P_\beta} Y_\beta$. Moreover, $S_\beta$ is coisotropic.
\end{lemma}
Suppose $z\in Z_\beta^{ss} = T^*W^{\beta(\Gm)}$.  
The infinitesimal $U^-$-action induces an injective map $\mathfrak{n}^-\rightarrow T_z(T^*W)$, and we get a direct sum $\mathfrak{n}^- \oplus T_z Y_{\beta} \subset T_z(T^*W)$.  Since $\beta$ acts on $U^-$ and hence compatibly on $\mathfrak{n}^-$, and $z$ is a $\beta$-fixed point, it makes sense to ask whether $\mathfrak{n}^-\rightarrow T_z(T^*W)$ is $\beta$-equivariant; it clearly is.  Hence the subspace $\mathfrak{n}^- \oplus T_z Y_{\beta}\subset T_z(T^*W)$ is $\beta$-invariant.   

Choose a $\beta$-invariant complementary subspace $V$, so 
$T_z(T^*W) = \mathfrak{n}^- \oplus T_z Y_{\beta} \oplus V$. 
  If $W$ is a $G$-representation, define 
  \bd
  N:= V\times Y_\beta = V\times Z_\beta^{ss}\times (T^*W)_+ \subset T^*W.
  \ed
  More generally, if $W$ is a smooth quasiprojective variety, let $W^\circ\subseteq W$ be a $\beta$-stable affine open subset containing $z$; one exists by \cite[Corollary~3.11]{Su}.  Further shrinking $W^\circ$ if necessary, let 
  \begin{equation}
  q: W^\circ\rightarrow \mathbb{A}^n
  \end{equation}
   be a $\beta$-equivariant \'etale map from a $\beta$-stable affine open subset $W^\circ\subseteq W$ with $q(z) = 0$, where $\mathbb{A}^n$ is a linear representation of $\Gm$; one exists by the \'Etale Slice Theorem (see p. 198 of \cite{Mumford}).  Since the map $q$ is \'etale, it induces a ``wrong way'' cotangent map  map $dq: T^*W^\circ \rightarrow T^*\mathbb{A}^n$ (note the slightly abusive notation). Defining $Z^\dagger_\beta = (T^*\mathbb{A}^n)^{\Gm}$ and $Y^\dagger_\beta\subset T^*\mathbb{A}^n$ to be the $\Gm$-attracting locus of $Z_\beta^\dagger$, then $Z_\beta\cap T^*W^\circ$ is a connected component of  $dq\inv(Z^\dagger_\beta)$ and $Y_\beta\cap T^*W^\circ$ is a connected component of $dq\inv(Y^\dagger_\beta)$.  The tangent map $d(dq_z): T_z(T^*W)\rightarrow T_0(T^*\mathbb{A}^n)$ is an isomorphism.  Abusively writing $V = d(dq_z)(V)\subset T_0(T^*\mathbb{A}^n) = T^*\mathbb{A}^n$, we get 
 \bd
 T^*\mathbb{A}^n = T^*_0\mathbb{A}^n = \mathfrak{n}^- \oplus V\oplus T_0 Y^\dagger_{\beta} = \mathfrak{n}^-\times V\times Y^\dagger_\beta.
 \ed
   Let $N$ denote the connected component of $dq\inv(V\times Y^\dagger_\beta) \subset T^*W^\circ\subseteq \bigvar$ containing $Y_\beta\cap T^*W^\circ$.  Note that $Y^\dagger_\beta$ is coisotropic, hence so is $V\times Y^\dagger_\beta$, hence so is $N$.

 \begin{prop}
 The subset $N$ above can be chosen to be a conical subset of $T^*W$.
 \end{prop}
 \begin{proof}
Writing $z\in T^*W$ as $z = (\pi(z),\xi)$ with $\pi(z) \in W$ and $\xi\in T_{\pi(z)}^*W$, the infinitesimal $U^-$-action $\mathsf{inf}(z): \mathfrak{n}^-\rightarrow T_z(T^*W) = T_{\pi(z)}W \oplus T_{\pi(z)}W^*$ decomposes as a sum of linear maps,  
\bd
\mathsf{inf}(z) = (b, f), \hspace{2em} b: \mathfrak{n}^- \rightarrow T_{\pi(z)}W, \hspace{2em} f: \mathfrak{n}^-\rightarrow T_{\pi(z)}^*W.
\ed
Note that $\beta$ has negative weights on $\mathfrak{n}^-$.  
With respect to the same decomposition, $T_zY_\beta \subset T_zT^*W$ is a direct sum of subspaces, $T_zY_\beta = T_{\pi(z)}W_{\geq 0} \oplus T_{\pi(z)}^*W_{\geq 0}$, where the subscript $\geq 0$ indicates the direct sum of subspaces of nonnegative $\beta$-weight.  Write $I = \on{Im}(b)\subseteq T_{\pi(z)}W$, a $\beta$-stable subspace, and choose a $\beta$-stable complement $C_1$ to $I\oplus T_{\pi(z)}W_{\geq 0}$ in $T_{\pi(z)}W$ (note that $I$ intersects the given subspace in $0$ since the subspace has nonnegative $\beta$-weights whereas the $\beta$-weights on $I$ are negative).  Furthermore, let
$K = \on{ker}(b)\subseteq \mathfrak{n}^-$, and choose a $\beta$-stable complement $C_2$ to $f(K)\oplus T_{\pi(z)}^*W_{\geq 0}$ in $T_{\pi(z)}^*W$ (again, the intersection of $f(K)$ with the given subspace must be zero). 

We claim that letting $V= C_1\oplus C_2$, that $T_zT^*W = \mathfrak{n}^-\oplus T_zY_\beta \oplus V$.  Note that 
\bd
\on{dim}(T_zT^*W) = \on{dim}(I) + \on{dim}(K) + \on{dim}(T_zY_\beta) + \on{dim}(V) = \on{dim}(\mathfrak{n}^-) + \on{dim}(T_zY_\beta) + \on{dim}(V)
\ed
by construction, so it suffices to show that $\mathfrak{n}^-\cap (V\oplus T_zY_\beta) = 0$.  Now $V\oplus T_zY_\beta = (C_1\oplus T_{\pi(z)}W_{\geq 0}) \oplus (C_2\oplus T_{\pi(z)}^*W_{\geq 0})$.  If the intersection is nonzero, then there is some element $n\in \mathfrak{n}^-$ with $(b(n),f(n))\in V\oplus T_zY_\beta$, which implies that $b(n)\in C_1\oplus T_{\pi(z)}W_{\geq 0}$ and $f(n)\in C_2\oplus T_{\pi(z)}^*W_{\geq 0}$.  The former condition implies that $b(n)=0$ by construction of $C_1$, so $n\in K$; thus $f(n) = 0$ as well.  This proves the claim.

Now, $V$ is by construction a conical subset of $T_zT^*W$.  The space $N$ constructed above then has the same property.
 \end{proof}

 Write $\on{pr}_\beta: Y_\beta\rightarrow Z_\beta$ for the projection. 
The following is Proposition 5.2 of \cite{McN}.  

\begin{prop}\label{existence of U_D}
For any $z \in Z_\beta^{ss}$ there are an affine neighborhood $D\subset Z_\beta^{ss}$ and a principal (conical) open subset $U_D\subset N$ (i.e. the complement of a hypersurface) such that
\begin{enumerate}
\item $U^-$ acts infinitesimally transversely to $U_D$; 
\item $\on{pr}_\beta\inv(D)\subseteq U_D \subseteq N$; and 
\item $(U^-\cdot U_D) \cap S_{>\beta} = \emptyset$. 
\item The complement $N\smallsetminus U_D$ is the hypersurface defined by a $\beta$-invariant function in $\C[N]$.
\item $U_D$ is coisotropic.
\end{enumerate}
The sets $D$ and $U_D$ can be chosen so that:
\begin{enumerate}
\item[(i)] $D$ and $U_D$ are $\beta$-stable. 
\end{enumerate}
 Furthermore, making, for each $z\in Z_\beta^{ss}$, a choice of any affine $D_{z}\subset Z_\beta^{ss}$ containing $z$  and any $U_{D_z}$ satisfying the conditions above,
 \begin{enumerate}
\item[(ii)] The union $\ds \bigcup_{z\in Z_\beta^{ss}} U^-\cdot U_{D_z}$  covers $U^-\cdot Y_\beta$. 
\end{enumerate}
\end{prop}
\begin{corollary}\label{cor:fflat}
Take any finite collection of points $z\in Z_\beta^{ss}$ and choices of $U_{D_z}$ as in the proposition.  Then the composite morphism of stacks
$\ds \coprod_{z\in Z_\beta^{ss}} U^-\times U_{D_z} \rightarrow \left(T^*X\smallsetminus S_{>\beta}\right) \rightarrow \left(T^*X\smallsetminus S_{>\beta}\right)/G$ is flat, and surjects onto a Zariski-open neighborhood of $S_\beta/G$.   In particular, the composite map is faithfully flat onto a neighborhood of $S_\beta/G$.
\end{corollary}
Since each map $U^-\times U_{D_z}\rightarrow T^*X$ is \'etale, it induces a symplectic structure on $U^-\times U_{D_z}$ for which the map is symplectic; that symplectic structure is described explicitly in \cite{McN}.  Moreover, one constructs \cite[Section~5]{McN} an \'etale and symplectic map 
\begin{equation}\label{eq:proj}
\phi_z: U^-\times U_{D_z}\longrightarrow T^*U^-\times \cS_z
\end{equation}
 to the product of $T^*U^-$ and a symplectic vector space $\cS_z$ with Hamiltonian $\beta$-action; the map $\phi_z$ is $\gK$-equivariant where $U^-$ acts trivially on $\cS_z$.

 \section{Quotients of DG Categories and Basic Equivalences}\label{sec:DGcatBasics}
 In this section we develop some basic properties of localization in dg categories.  A more comprehensive and broader treatment of parts of this section and the next section can be found in the appendix of \cite{BNP}.
 \subsection{Quotients of DG Categories}\label{sec:quotients}
 The following assumptions on a dg category $D$ and full subcategory $K$ are in force for the remainder of Sections \ref{sec:DGcatBasics} and \ref{sec:secondDG}.
 \begin{assumptions}\label{ass:1}
 Let $D$ be a pre-triangulated dg category which is locally presentable with all small (homotopy) colimits.  Let $K\subset D$ be a full dg subcategory, closed under shifts, and write $i_*: K\hookrightarrow D$ for the inclusion.
 \end{assumptions}
 Following Keller, Drinfeld defined a dg quotient category $D/K$ and projection functor, which we will denote $j^*: D\rightarrow D/K$ \cite{Drinfeld}.

 The following must be standard, but we do not know a reference.
 \begin{prop}\label{prop:quotient-properties}
 \mbox{}
 Suppose that, in addition to the properties above, $K$ is closed under taking (small homotopy) colimits in $D$.  Then:
 \begin{enumerate}
 \item A right adjoint $i^!$ to the inclusion $i_*: K\rightarrow D$ exists, and $i^!i_*\simeq \on{Id}_K$.
 \item A right adjoint $j_*$ to the projection $j^*$ exists.
 \item  $j^*j_*  \simeq \on{Id}_{D/K}$; in particular, $j_*$ is faithful.
  \item One gets an exact triangle of functors, $i_*i^! \rightarrow \on{Id} \rightarrow j_*j^* \xrightarrow{[1]}$.
 \end{enumerate}
 \end{prop}
 \begin{proof}
 Since $D$ contains all colimits and $K$ is full and closed under colimits in $D$, the Adjoint Functor Theorem implies that $i_*$ has an exact right adjoint $i^!$.  By fullness of $K$, $i^!i_*\simeq \on{Id}_K$.  This proves (1).
 
For an object $M$ of $D$, consider the cone $c(M) = \on{cone}(i_*i^!M\rightarrow M)$; this defines a functor on $D$.  By the previous
paragraph, $c\circ i_*\simeq 0$, and thus, by the universal property of the quotient $D/K$, the functor $c: D\rightarrow D$ factors through the projection $j^*: D\rightarrow D/K$.  Write $c = j_*j^*$.  Assertion (4) is immediate.  Applying $j^*$ to the exact triangle
$i_*i^! \rightarrow \on{Id}\rightarrow j_*j^*\xrightarrow{[1]}$ and using that $j^*i_*i^!\simeq 0$ we get $j^*j_*j^*\simeq j^*$; since $j^*$ is full and essentially surjective, (3) follows.  

Because $i^!i_*\simeq \on{Id}_K$, any morphism
$i_*M\rightarrow N$ factors through $i_*i^!N\rightarrow N$ via the unit of adjunction.  Because the functor $c$ is the cone on the unit of adjunction $i_*i^!\rightarrow\on{Id}_D$, we find that $\Hom(i_*M',j_*j^*N)\simeq 0$ for all $M'\in \on{ob}(K)$ and $N\in \on{ob}(D)$.  

We next prove that $j_*$ is right adjoint to $j^*$. 
Consider the maps
\bd
\Hom(M,j_*j^*N)\xrightarrow{j^*} \Hom(j^*M, j^*N)\xrightarrow{j_*} \Hom(j_*j^*M, j_*j^*N)\xrightarrow{j^*} \Hom(j^*M, j^*N).
\ed
Since $j^*j_*$ is the identity, to prove that $j_*$ is right adjoint to $j^*$ it suffices to prove that the left-hand composite above is an isomorphism.  Using the exact triangle $i_*i^!M \rightarrow M \rightarrow j_*j^*M\xrightarrow{[1]}$ in the first variable, this will follow if $\Hom(i_*i^!M, j_*j^*N)\simeq 0$, which is immediate from the previous paragraph (using $M'=i^!M$).
\end{proof}
 
We write $j_!$ for the left adjoint of $j^*$: this is a partially defined functor, and we wish to know when it is defined.
    \begin{lemma}\label{lem:extension}
Let $M$ be an object of $D/K$ and suppose that $\overline{M}$ is any extension to an object of $D$: that is, 
$j^*\overline{M} \simeq M$.  If $\Hom(\overline{M}, N)=0$ for all $N$ for which $j^*N = 0$, then $j_!M$ is defined and $j_!M\cong \overline{M}$.  Moreover, such an assignment $M\mapsto \overline{M}$ is functorial: given $\overline{M}_1$ and $\overline{M}_2$ satisfying $\Hom(\overline{M}_i, N)=0$ for all $N$ for which $j^*N = 0$, then the natural map 
\bd
\Hom(\overline{M}_1, \overline{M}_2)\rightarrow \Hom(j^*\overline{M}_1, j^*\overline{M}_2)
\ed
is an isomorphism.
\end{lemma}
\begin{proof}
We use the exact triangle $i_*i^! N\rightarrow N \rightarrow j_*j^*N \xrightarrow{[1]}$.  Applying $\Hom(\overline{M},-)$, we get
\bd
0=\Hom(\overline{M},i_*i^!N)\rightarrow \Hom(\overline{M},N) \rightarrow \Hom(\overline{M},j_*j^*N)\xrightarrow{[1]},
\ed
where $\Hom(\overline{M},i_*i^!N)$ vanishes since $j^*\circ i_* = 0$.  Thus, using the adjunction 
\bd
\Hom(\overline{M}, j_*j^*N) = \Hom(j^*\overline{M},j^*N) = \Hom(M, j^*N), 
\ed
 we get $\Hom(\overline{M},N) = \Hom(M,j^*N)$, proving the first assertion.  For the second assertion, use the exact triangle $i_*i^!\overline{M}_2\rightarrow \overline{M}_2\rightarrow j_*j^*\overline{M}_2$ in the second argument of $\Hom$.  
\end{proof}

 \begin{prop}\label{prop:iupperstar}
 Suppose that the left adjoint functor $j_!$ to $j^*$ exists (i.e., is fully defined).  Then $i_*$ has a left adjoint $i^*$ sitting in an exact triangle
 $ j_!j^*\rightarrow \on{Id}\rightarrow i_*i^*\xrightarrow{[1]}.$
 \end{prop}
 \begin{proof}
 By Lemma \ref{lem:extension} we have $j^*j_!\simeq \on{Id}$.  It follows that the cone of $j_!j^*\rightarrow\on{Id}$ factors through $i_*$; write $i_*i^*$ for the factorization.  Using the exact triangle in the first variable of $\Hom$ and the fact that the image of $j_!$ is left orthogonal to $K$ gives
 \bd
 \Hom(i_*i^*M, i_*N)\xrightarrow{\simeq} \Hom(M,i_*N).
 \ed
 Since $i_*:K\rightarrow D$ is a full inclusion, this yields the adjunction.    
 \end{proof}
  
 \subsection{Homotopy Limits}\label{sec:homotopy-limits}
 We next study co-simplicial diagrams of categories.
 \begin{assumptions}\label{ass:cosimplicial}
 Suppose that $\{D[n]\}$ is semi-cosimplicial diagram of dg categories, i.e., a diagram as follows:
 \begin{equation}\label{eq:semicosimplicial}
\xymatrix{D[0] \ar@<0.5ex>[r]\ar@<-0.5ex>[r] &  D[1] \ar@<0.7ex>[r]\ar@<-0.7ex>[r]\ar[r]
& D[2] \ar@<0.9ex>[r] \ar@<0.3ex>[r] \ar@<-0.3ex>[r] \ar@<-0.9ex>[r] & \dots},
\end{equation}
satisfying standard axioms.
We assume that each $D[n]$ satisfies Assumptions \ref{ass:1}, and that all functors are continuous (colimit-preserving).  Suppose that each $D[n]$ is equipped with a full dg subcategory $K[n]$ closed under colimits in $D[n]$; write $D[n]/K[n]$ for the DG quotient.  Suppose that the structure functors in the semi-cosimplicial diagram take each $K[n]$ into $K[n+1]$.  
\end{assumptions}
It follows that each inclusion $i_{[n]\ast}: K[n]\rightarrow D[n]$ and quotient $j_{[n]}^*: D[n]\rightarrow D[n]/K[n]$ has a right adjoint, $i_{[n]}^!$, respectively $j_{[n]\ast}$.  
\begin{lemma}\label{holims of quotients}
Let $D$ denote the limit of $\{D[n]\}$, $K$ denote the limit of $\{K[n]\}$, and $D/K$ denote the limit of $\{D[n]/K[n]\}$.  Then:
\begin{enumerate}
\item $D$ has all small colimits.
\item The natural sequence $K\rightarrow D\rightarrow D/K$ is a homotopy fiber sequence.
\item The natural functors $\pi^!: D\rightarrow D[0]$, $\pi^!: K\rightarrow K[0]$, $\pi^!: D/K\rightarrow D[0]/K[0]$ are faithful.
\item $K$ is a full subcategory of $D$, closed under colimits in $D$.  
\end{enumerate}
\end{lemma}
It follows from the lemma and Proposition \ref{prop:quotient-properties} that the inclusion $i_*:K\rightarrow D$ and the quotient
$j^*: D\rightarrow D/K$ have right adjoints $i^!$ and $j_*$, respectively.
\begin{prop}\label{prop:descent-of-adjoints}
\mbox{}
\begin{enumerate}
\item The functors $i^!_{[n]}: D[n]\rightarrow K[n]$ and $j_{[n]\ast}: D[n]/K[n]\rightarrow D[n]$ assemble into functors of semi-cosimplicial diagrams of DG categories.
\item Define a functor $\wt{i}^!: D\rightarrow K$ as the homotopy limit of the compatible system $i^!_{[n]}$.  Then $\wt{i}^!\simeq i^!$.
\item Define a functor $\wt{j}_*: D/K\rightarrow D$ as the homotopy limit of the compatible system $j_{[n]\ast}$.  Then $\wt{j}_*\simeq j_*$.  
\end{enumerate}
\end{prop}
 \begin{corollary}
 We have $\pi^!j_*\simeq j_*\pi^!$.
 \end{corollary}

 \subsection{Basic Equivalences}
 We maintain the set-up of Section \ref{sec:quotients}.  We suppose, in addition:
 \begin{assumption}\label{ass:cg}
 Assume that $D$ is compactly generated. 
\end{assumption} 
We suppose $K\subset D$ is a full subcategory closed under colimits in $D$.  When we want to emphasize the inclusion we write $i_*: K\hookrightarrow D$ for the functor.  
Thus we have a diagram:
\bd
\xymatrix{
K  \ar@/^0.5pc/[r]^{i_*}  & 
D \ar@/^0.5pc/[l]^{i^!}
%\ar@/_1pc/[l]_{i^*} 
\ar@/^0.5pc/[r]^{j^*}  
& D/K \ar@/^0.5pc/[l]^{j_*}.
%\ar@/_1pc/[l]_{j_!}
}
\ed
 
 As in Section \ref{sec:quotients}, we have an exact triangle of functors
 \begin{equation}\label{basic exact triangle}
 i_*i^! \longrightarrow \on{Id}_D \longrightarrow j_*j^* \xrightarrow{[1]}.
 \end{equation}

 \begin{defn} A collection $S$ of objects of a dg category $\cC$ {\em Karoubi-generates} $\cC$ if every object in the homotopy category of $\cC$
can be
obtained from objects in $S$ by
finite iteration of operations of taking the cone of a morphism,
and passing to a direct summand of an object.
\end{defn}
\begin{notation}
For a dg category $\cC$, we write $\cC^c$ for the full subcategory of compact objects of $\cC$.
\end{notation}
\begin{lemma}\label{lem:karoubigenerate}
Let $j^*: D\rightarrow D/K$ be the quotient functor.  Assume that the right adjoint $j_*$ is continuous.  Then:
\begin{enumerate}
\item For every $M\in D^c$, one has $j^*M\in (D/K)^c$.  
\item The subcategory $\{j^*M \, |\, M\in D^c\}$ Karoubi-generates $(D/K)^c$.
\end{enumerate}
\end{lemma}
\begin{proof}
The functor $j^*$ has a continuous right adjoint by assumption and Proposition \ref{prop:quotient-properties}, hence (1) follows by \cite{Neeman} or \cite[Proposition~1.2.4]{DG}. Since $j^*$ is also faithful, any collection of compact generators Karoubi-generates the category $(D/K)^c$ of compact objects in $D/K$ by \cite[Corollary~1.4.6]{DG}.  
\end{proof}

\begin{prop}\label{basic equivalences}
Assume that $j_*$ is continuous.  The following are equivalent.
\begin{enumerate}
\item[(i)] The functor $i^!$ takes compact objects to compact objects.
\item[(ii)] The functor $j_*$ takes compact objects to compact objects.
\item[(ii$^\prime$)] The functor $j_*$ takes a collection of compact generators to compact objects. 
\item[(iv)] The functor $i^!$ admits a continuous right adjoint.
\end{enumerate}
Moreoever, the following conditions are equivalent to each other:
\begin{enumerate}
\item[(iii)] The functor $j_!$ left adjoint to $j^*$ is defined on all of $D/K$.
\item[(iii$^\prime$)] The functor $j_!$ left adjoint to $j^*$ is defined on all compact objects of $D/K$.
\item[(iii$^{\prime\prime}$)] The functor $j_!$ left adjoint to $j^*$ is defined on a set of compact generators of $D/K$. 
\end{enumerate}
\end{prop}
The proposition is a direct analogue of Proposition 2.1.2 of \cite{DG}.  We include a proof for completeness, and since it differs in some places from \cite{DG}.  Later (Proposition \ref{prop:full-equiv}) we will establish an equivalence of conditions (iii)-(iii$^{\prime\prime}$) with the other conditions using some additional duality assumptions.
\begin{proof}
 If (ii) holds, then for any compact $\cF\in D$, $j_*j^*\cF$ is also compact by Lemma \ref{lem:karoubigenerate}, hence by \eqref{basic exact triangle} so is $i_*i^!\cF$ (cones of morphisms of compact objects are compact).  Also   
 $i_*$ is fully faithful, hence $i^!\cF$ is compact as well.  Thus $\text{(ii)}\implies\text{(i)}$.  Note that by definition/construction the functor $j^*$ is essentially surjective.  By Section 0.8.8 of \cite{DG}, it follows that the image under $j^*$ of the category of compact objects $D^c$ of $D$ Karoubi generates the category of compact objects $(D/K)^c$ of $D/K$.  Now suppose $\cF$ is compact in $D/K$; by the previous sentence, we may write $\cF$ as a direct summand of some $j^*\cF'$ for a compact $\cF'$.
   If (i) holds, then $i_*i^!\cF'$ is compact, hence so is $j_*j^*\cF'$; since $j_*\cF$ is a direct summand of a compact object, it is thus compact, proving (ii).  
 
Using Lemma \ref{lem:karoubigenerate}, we see that (ii$^\prime$), and the facts that cones of morphisms of compact objects and direct summands of compact objects are compact, implies (ii).  (ii)$\implies$(ii$^\prime$) is clear.
 
Since $i_*$ preserves compactness, $i^!$ is continuous by \cite{Neeman}.  Thus (i)$\iff$(iv) by \cite{Neeman}.  

We have (iii)$\iff$(iii$^\prime$)$\iff$(iii$^{\prime\prime}$) since $D/K$ is compactly generated: supposing (iii$^{\prime\prime}$) holds, for any colimit over a collection of compact generators of $D/K$, define $j_!$ by continuity.  
 \end{proof}

\begin{prop}\label{prop:shriekdefined}
Suppose that, for every object $M$ in a collection of compact generators of $D/K$, there exists a compact object 
$\overline{M}$ of $D$ such that $\Hom(\overline{M}, N)=0$ for all $N$ for which $j^*N=0$.  Then $j_!$ is defined on all of
$D/K$.
\end{prop}
\begin{proof}
By Lemma \ref{lem:extension}, the hypotheses imply that $j_!M$ is defined for all $M$ in a collection of compact generators.  Hence by
Proposition \ref{basic equivalences}, $j_!$ is defined on all objects.
\end{proof} 

Now suppose $K'\subseteq K\subseteq D$ are full subcategories closed under colimits in $D$.  Let $\overline{j}^*: D\rightarrow D/K$, $j^*: D/K'\rightarrow D/K$ and $\iota^*: D\rightarrow D/K'$ be the quotient functors.  For a quotient functor $p$, we use the notation $p_!$, $p_*$ for the left and right adjoints to $p^*$, where $p_!$ may be only partially defined.   
\begin{lemma}\label{lem:jbarimpliesj}
Suppose that $M\in D/K$ and that $\overline{j}_!M$ exists.  Then $j_!M$ exists, and $j_!M = \iota^*\overline{j}_!M$.
\end{lemma}
\begin{proof}
By Lemma \ref{lem:extension}, we have $\overline{j}^*\overline{j}_!M = M$.  Hence $j^*(\iota^*\overline{j}_!M) = M$.  Thus, again by Lemma \ref{lem:extension}, it suffices to show that $\Hom(\iota^*\overline{j}_!M,N) = 0$ for all $N\in K/K'$.  Now
\bd
\Hom(\iota^*\overline{j}_!M,N) = \Hom(\overline{j}_!M,\iota_*N) = \Hom(M, \overline{j}^*\iota_*N) = \Hom(M, j^*\iota^*\iota_*N) = \Hom(M,j^*N)=0,
\ed
where the last equality holds because $K/K'$ is the kernel of $j^*$.  
\end{proof}
\begin{prop}\label{prop:presofcompactness}
If $\overline{j}_*$ preserves compactness, then $j_*$ preserves compactness.
\end{prop}
\begin{proof}
If $\overline{j}_*$ preserves compactness, then by Proposition \ref{basic equivalences}, $\overline{j}_!$ is defined for all objects.  Hence by Lemma
\ref{lem:jbarimpliesj}, $j_!$ is defined for all objects.  Another application of Proposition \ref{basic equivalences} implies that $j_*$ preserves compactness.
Proposition \ref{prop:presofcompactness} follows.
\end{proof}

\section{$t$-Structures, Coherent Objects, and Verdier Duality}\label{sec:secondDG}
We maintain the setting of Section \ref{sec:DGcatBasics}; in particular, Assumptions \ref{ass:1}, \ref{ass:cosimplicial}, \ref{ass:cg} are in force whenever they make sense.  The main purpose of this section is to show that, in the presence of a $t$-structure with a reasonable heart, a Verdier-duality--type functor, and a projection formula, the conditions of Proposition \ref{basic equivalences} can be established by ``preservation of coherence.''  We will use this later to check preservation of coherence for twisted $\D$-modules by reduction to a model case as in \cite{McN}.
 \subsection{$t$-Structures and Abelian Categories}
 We continue with a pre-triangulated dg category $D$.
 \begin{assumption}\label{ass:t-str}
 Suppose $D$ comes equipped with a fixed nondegenerate $t$-structure.
 \end{assumption}
   If $D$ is a dg category 
 with $t$-structure, we write $D^\heartsuit$ for the abelian category that is the heart of the $t$-structure.  We write 
 $H^0: D\rightarrow D^\heartsuit$ for the cohomology functor and assume that it commutes with colimits.
 
 Suppose $\{D[n]\}$ is a cosimplicial diagram of dg categories as in Section \ref{sec:homotopy-limits}.  
 \begin{assumptions}\label{ass:hearts}
 Each heart $D[n]^\heartsuit$ is assumed to be a locally Noetherian Grothendieck abelian category.  We assume that the structure functors of the semi-cosimplicial diagram \eqref{eq:semicosimplicial} are 
 $t$-exact and quasi-faithful.

 We further assume that the hearts $D[n]^\heartsuit$ come equipped with a compatible collection of localizing subcategories, $D[n]^\heartsuit_S$: by ``compatible,'' we mean that the induced functors $D[n]^\heartsuit\rightarrow D[n+1]^\heartsuit$ take $D[n]^\heartsuit_S$ into $D[n+1]^\heartsuit_S$ for all $n$ (and each structural functor).  In particular, each such subcategory of the heart is closed under colimits.
  \end{assumptions}
 We then equip the semi-cosimplicial diagram $\{D[n]\}$ with a semi-cosimplicial diagram $\{K[n]\}$ of full subcategories, closed under colimits in $D[n]$, each defined as the full subcategory of objects whose cohomologies lie in the corresponding $D[n]^\heartsuit_S$.  The quotient $D[n]/K[n]$ then inherits a $t$-structure making $j^*_{[n]}$ $t$-exact.
 \begin{prop}
 \mbox{}
 \begin{enumerate}
 \item For each $n$, one obtains an exact sequence 
 \bd
 0\rightarrow K[n]^\heartsuit \rightarrow D[n]^\heartsuit \rightarrow (D[n]/K[n])^\heartsuit \rightarrow 0
 \ed
 of abelian categories.
 \item The homotopy limits $K$, $D$, and $D/K$ inherit $t$-structures for which the functors $\pi^!: K\rightarrow K[0]$, $\pi^!: D\rightarrow D[0]$, $\pi^!: D/K\rightarrow D[0]/K[0]$ are $t$-exact.
 \item The hearts $K^\heartsuit$, $D^\heartsuit$, $(D/K)^\heartsuit$ are naturally isomorphic to the  limits of the hearts $K[n]^\heartsuit$, $D[n]^\heartsuit$, $(D[n]/K[n])^\heartsuit$ respectively, compatibly with the cohomology functors.
 \end{enumerate}
 \end{prop}
 \begin{defn}
 An object of $D$ is {\em coherent} if it has only finitely many nonzero cohomology objects in $D^\heartsuit$, and each of those cohomology objects is a Noetherian object of the abelian category $D^\heartsuit$.
 \end{defn}
 The full subcategory of $D$ consisting of coherent objects is denoted $D_{\on{coh}}$.
 \begin{lemma}\label{lem:reflect-coherence}
 Under the above conditions, for $M\in D$, if $\pi^!M\in D[0]_{\on{coh}}$ then $M\in D_{\on{coh}}$.  
 \end{lemma}
 
 \begin{prop}\label{lem:cptimpliescoh}
Suppose that $\pi^!: D\rightarrow D[0]$ has a continuous right adjoint $\pi_*$ and that if an object $N \in D[0]$ is compact then it is coherent.  If $M\in D^c$, then $M\in D_{\on{coh}}$.  
\end{prop}
\begin{proof}
Since $\pi_*$ is continuous, $\pi^!$ preserves compactness.  Thus, if $M\in D$ is compact, so is $\pi^!M$; hence by assumption $\pi^!M$ is coherent.  The proposition then follows by Lemma \ref{lem:reflect-coherence}.
\end{proof}
\begin{assumption}\label{ass:c-implies-coh}
For the remainder of the section, we assume that 
\bd
D^c\subseteq D_{\on{coh}}\hspace{2em} \text{and} \hspace{2em} (D/K)^c\subseteq (D/K)_{\on{coh}}.
\ed
\end{assumption}

\subsection{Verdier Duality}\label{sec:formalVD}
We will want some formal properties of Verdier duality.  
Thus, we assume given two categories $D$, $D'$, each equipped with all the additional structures and properties from 
Assumptions \ref{ass:1}, \ref{ass:cosimplicial}, \ref{ass:cg}, \ref{ass:t-str}, and \ref{ass:hearts}.  In particular, each of $D$, $D'$ is the colimit of a cosimplicial diagram.  We assume also given
an exact equivalence $\vD: (D'_{\on{coh}})^{\on{op}}\rightarrow D_{\on{coh}}$, which we refer to as a {\em Verdier duality functor}, coming from a compatible system of such equivalences on the terms of the cosimplicial diagrams defining $D$ and $D'$.  We usually abusively write $\vD$ also for the quasi-inverse functor.  
\begin{remark}
Writing $\pi^!: D\rightarrow D[0]$ and also $\pi^!: D'\rightarrow D'[0]$, we have
$\vD\circ\pi^!\simeq \pi^!\circ\vD$.
\end{remark}

\subsection{Tensor Product}
We next want to establish some formal properties of Verdier duality and tensor products.  

Thus, we assume given two categories $D$, $D'$, each equipped with all the additional structures and properties of
Assumptions \ref{ass:1}, \ref{ass:cosimplicial}, \ref{ass:cg}, \ref{ass:t-str}, and \ref{ass:hearts} and Section \ref{sec:formalVD}.  In particular, each of $D$, $D'$ is the colimit of a cosimplicial diagram.    

Let $\on{Vect}_{\C}$ denote the dg category of complex dg vector spaces.  We define a bi-functor 
\bd
D'_{\on{coh}}\times D \rightarrow \on{Vect}_{\C} \hspace{4em} \text{by}
\ed
\bd
\ds (M, N)\mapsto M\underset{D}{\otimes} N := \Hom_{D}(\vD(M),N).
\ed
Recall that we are assuming given subcategories $K'\subseteq D'$, $K\subseteq D$ as before.  Suppose that, in addition:
\begin{assumption}\label{ass:cohpres}
$\vD(K'_{\on{coh}}) =  K_{\on{coh}}$.  
\end{assumption}
Then the above functor induces a compatible functor $(D'/K')_{\on{coh}}\times D/K \rightarrow \on{Vect}_{\C}$ defined by
\bd
\ds (M, N)\mapsto M\underset{D/K}{\otimes} N := \Hom_{D/K}(\vD(M),N).
\ed
\begin{prop}
Fixing $M\in (D'/K')_{\on{coh}}$ and $N\in D/K$, we have that $M\underset{D/K}{\otimes} N$ is the homotopy limit of either row of 
\bd
\xymatrix{\pi^!M\underset{D[0]/K[0]}{\otimes}\pi^!N \ar[d]_{\simeq} \ar@<0.5ex>[r]^{a_c^!}\ar@<-0.5ex>[r]_{p^!}  & p^!\pi^!M\underset{D[1]/K[1]}{\otimes} p^!\pi^!N \ar[d]^{\simeq} \ar@<0.7ex>[r]\ar@<-0.7ex>[r]\ar[r] & \dots\\
 \Hom_{D[0]/K[0]}(\vD(\pi^!M), \pi^!N) \ar@<0.5ex>[r]^{a_c^!}\ar@<-0.5ex>[r]_{p^!}  &  \Hom_{D[1]/K[1]}(\vD(p^!\pi^!M), p^!\pi^!N) \ar@<0.7ex>[r]\ar@<-0.7ex>[r]\ar[r] & \dots }
\ed
\end{prop}
 
\begin{prop}\label{coherence is enough}
Assume the following:
\begin{enumerate}
\item $\pi^!$ has continuous right adjoint $\pi_*$.  
\item A coherent object $M$ is compact if and only if $\vD(M)$ is compact.
\item For $M\in (D'/K')_{\on{coh}}$, we are given a quasi-isomorphism of functors of $N$, 
\bd
j_*j^*M\otimes_{D/K} N \simeq M\otimes_{D/K} j_*j^*N.
\ed
\item $j_*: D'/K' \rightarrow D'$  is continuous and preserves coherence: that is, 
\bd
\text{if}\,\,\, M \in (D'/K')_{\on{coh}} \,\, \text{then}\,\, j_*M \in D'_{\on{coh}}.
\ed
\end{enumerate}
Then $j_*: D'/K' \rightarrow D'$ preserves compactness as well.  
\end{prop}
\begin{proof}
First, consider a compact object $j^*M$ for some compact object $M$ of $D'/K'$.  By assumption (4) and Proposition \ref{lem:cptimpliescoh}, $j_*j^*M$ is coherent.  Hence 
\begin{equation}\label{eq:continuity}
\Hom_{D/K}(\vD(j_*j^*M), N)= j_*j^*M\otimes_{D} N \overset{(\dagger)}{=} M\otimes_{D}j_*j^*N = \Hom_{D}(\vD(M), j_*j^*N),
\end{equation}
where $(\dagger)$ follows from assumption (3).
Now $j^*$ and $j_*$ are continuous by construction and assumption (4), and $\vD(M)$ is compact by assumption (2), so the final functor of \eqref{eq:continuity} is continuous.  Hence
$\vD(j_*j^*M)$ is compact, which by assumption (2) implies that $j_*j^*M$ is compact.  Since objects $j^*M$ for compact $M$ Karoubi-generate 
$(D'/K')^c$ by Lemma \ref{lem:karoubigenerate}, the corollary follows.
\end{proof}

 \begin{prop}\label{prop:cpt->coh}
 Assume that:
 \begin{enumerate}
 \item For every object $M\in (D/K)^\heartsuit \cap (D/K)_{\on{coh}}$ there exist an object $M_0\in (D/K)^{\leq 0} \cap (D/K)^c$ and a morphism $M_0\rightarrow M$ that induces a surjection $H^0(M_0)\rightarrow H^0(M)$.  
 \item The functor $j_*$ has finite cohomological amplitude.
 \item The functor $j_*$ is left exact.
 \item The functor $j_*: D/K \rightarrow D$ preserves compactness. 
 \end{enumerate} Then $j_*$ preserves coherence.
 \end{prop}
 \begin{proof} 
 The proof follows closely that of \cite[Proposition~3.5.2]{DG}. As in Lemma 9.4.7 of \cite{DG2}, it follows from the first hypothesis that for every $k\in\mathbb{Z}$ and every $M\in (D/K)_{\on{coh}}$ there are an $M_k\in (D/K)^c$ and a morphism $\on{app}_k: M_k\rightarrow M$ whose cone belongs to $(D/K)^{\leq -k}$.  
 
 By the second hypothesis, $j_*$ has a finite cohomological amplitude $A$.  For $M\in (D/K)^{[a,b]}\cap (D/K)_{\on{coh}}$, choose $M_{a-A-1}$ as in the previous paragraph.  We get an exact triangle 
 \bd
 j_* M_{a-A-1} \rightarrow j_*M \rightarrow j_*\on{Cone}(\on{app}_{a-A-1}) \xrightarrow{[1]}.
 \ed
 By assumptions (2) and (3), $j_*M \in D^{[a,b+A]}$ whereas $j_*\on{Cone}(\on{app}_{a-A-1})$ lies in $D^{\leq a-1}$.  Thus 
 $j_* M_{a-A-1} \rightarrow j_*M$ induces an isomorphism on $H^k$ for every $k\in [a,b+A]$.  Since $j_*M_{a-A-1}$ is compact by assumption (4), its cohomologies are coherent (by our global assumption on the compact subcategory of $D$), proving the assertion.
 \end{proof}

\begin{prop}\label{prop:full-equiv}
Assume given $K\subseteq D$, $K'\subseteq D'$, $\vD: D_{\on{coh}}\rightarrow D'_{\on{coh}}$ as above.  Suppose that 
$j_*: D/K\rightarrow D$ is continuous, and that 
$j_*: D'/K'\rightarrow D'$
is continuous and satisfies conditions (1), (2), and (3) of Proposition \ref{prop:cpt->coh} (where $D$ is replaced by $D'$ in that proposition and $K$ is replaced by $K'$).  
Suppose also that the hypotheses of Proposition \ref{lem:cptimpliescoh} hold for both $D$ and $D'$.  
Then all conditions (i) through (iv) of Proposition \ref{basic equivalences} for $D$ are equivalent.
\end{prop}
\begin{proof}
Assume Proposition \ref{basic equivalences}(ii) holds.  Let $M$ be a compact object of $D/K$.  Then $j_*\vD(M)$ is coherent by Proposition \ref{prop:cpt->coh}.
Hence $\vD\big(j_*\vD(M)\big)$ is both compact and coherent.  Now for any coherent object $N$ of $D$, 
\bd
\begin{split}
\Hom_{D}\big(\vD\big(j_*(\vD(M)\big), N\big) = \Hom_{D'}\big(\vD(N), j_*\vD(M)\big) 
  = \Hom_{D'/K'}\big(j^*\vD(N), \vD(M)\big) = \\
  \Hom_{D/K}\big(M, \vD(j^*\vD(N))\big) 
 = \Hom_{D/K}(M, j^*N),
 \end{split}
\ed
where the last equality follows since $j^*$ commutes with duality.  Since $j^*$ is continuous and both $M$ and $\vD\big(j_*\vD(M)\big)$ are compact, both sides commute with colimits in $N$.  Hence $j_!M = \vD\big(j_*\vD(M)\big)$ is defined, proving Proposition \ref{basic equivalences}(iii$^\prime$).

Conversely, suppose Proposition \ref{basic equivalences}(iii$^\prime$) holds, i.e., $j_!$ is defined for compact objects of $D/K$.  Since $j_!$ is left adjoint to a continuous functor, it preserves compactness; hence for compact $M$ in $D/K$, $\vD(j_! \vD(M))$ is also compact.  Given a compact object $M$ of $D/K$ and a coherent object $N$ of $D$, then 
\bd
\begin{split}
\Hom\big(N, \vD(j_!\vD(M))\big) = \Hom\big(j_!\vD(M), \vD(N)\big) 
= \Hom\big(\vD(M), j^*\vD(N)\big) \\
 = \Hom\big(\vD(M), \vD(j^*N)\big)
= \Hom(j^*N,M).
\end{split}
\ed
Since $j^*$ is continuous and every $N'$ is a colimit of compact, hence coherent, objects $N$, we find that $j_*M = \vD(j_!\vD(M))$.
Since the latter is compact, we conclude that Proposition \ref{basic equivalences}(ii) holds.  Combined with the other equivalences of 
Proposition \ref{basic equivalences}, this completes the proof.
\end{proof}
\begin{remark}
We remark that when the equivalent conditions of Proposition \ref{basic equivalences} hold, the proof of Proposition \ref{prop:full-equiv} shows that $j_! = \vD\circ j_*\circ \vD$; in particular, $j_*$ preserves coherence if and only if $j_!$ does.
\end{remark}

\section{Microlocal Categories and Basic Functors}
\subsection{Twisted $\D$-Modules on Stacks}
If $X$ is a smooth variety (smooth $\C$-scheme of finite type), we let $\D(X)$ denote the unbounded, quasicoherent derived category of $\D$-modules on $X$.  
Let $X/G$ be a smooth stack, the quotient of a smooth quasi-projective variety $X$ by an affine algebraic group $G$.  Let $\pi: X\rightarrow X/G$ denote the projection and let $\mu: \mathfrak{g}\rightarrow T_X$ denote the infinitesimal $G$-action.  
 Let $c: \mathfrak{g}\rightarrow \C$ denote a Lie algebra character.  
 
Associated to the action and projection maps $a, p: G\times X\rightarrow X$ there are pullback functors $a^!, p^!: \D(X)\rightrightarrows \D(G\times X)$; extending further, we get the standard co-simplicial diagram for descent of $\D$-modules.  We twist this diagram, modifying only the right-most (action) maps, using $c$, as follows.  For the action map $a:G\times X\rightarrow X$, consider the usual bimodule $\D_{G\times X\rightarrow X}$ with left action by
$\D_{G\times X} = \theo_G\otimes U(\mathfrak{g})\otimes \D_X$; we twist the left module structure by letting elements of $\mathfrak{g}$ act via $\mu+c$ rather than $\mu$.  This new bimodule structure yields a new pullback functor which we denote by $a_c^!$.  Similarly, we may replace the usual left action of $\D(G^{n-1}\times G\times X)$ on the bimodule $\D_{G^n\times X \rightarrow G^{n-1} \times X}$ associated to the map
$G^{n-1}\times G\times X \xrightarrow{1_{G^{n-1}}\times a} G^{n-1}\times X$ by letting vector fields on the right-most copy of $G$ act via $\mu+c$.  
Thus building from $a_c^!$ and $p^!$, we obtain a semi-cosimplicial diagram \eqref{twisted-diagram}. 

\begin{defn}\label{def:D-mod}
The category $\D(X/G,c)$ is the homotopy limit of the diagram
\begin{equation}\label{twisted-diagram}
\xymatrix{\D(X) \ar@<0.5ex>[r]^{\hspace{-1em}a_c^!}\ar@<-0.5ex>[r]_{\hspace{-1em}p^!} &  \D(G\times X) \ar@<0.7ex>[r]\ar@<-0.7ex>[r]\ar[r]
& \D(G\times G \times X) \ar@<0.9ex>[r] \ar@<0.3ex>[r] \ar@<-0.3ex>[r] \ar@<-0.9ex>[r] & \dots}
\end{equation}
 of dg categories.  We write $\D(X/G,c)\xrightarrow{\pi^!} \D(X)$ for the natural functor.
\end{defn}
\begin{remark}
When the moment map $\mu:T^*X\rightarrow\mathfrak{g}^*$ is flat, then, as in Proposition 4.10 of \cite{McNderived}, the category $\D(X/G,c)$ is canonically equivalent to the canonical dg enhancement (cf. \cite{Keller, Drinfeld, ToMorita}) of the unbounded derived category of the abelian category $\on{Qcoh}(\D,G,c)$ of $(G,c)$-equivariant quasicoherent $\D$-modules.  

In general, for each $i$ there is an $i$th cohomology functor $$\mathcal{H}^i: \D(X/G,c)\rightarrow \on{Qcoh}(\D,G,c).$$
\end{remark}
\begin{lemma}\label{lem:pi-functors}
The functor $\D(X/G,c)\xrightarrow{\pi^!} \D(X)$ has continuous right adjoint $\pi_*: \D(X)\rightarrow\D(X/G,c)$.  
\end{lemma}
\begin{proof}
In the untwisted case, \cite[Section~6.1.7]{DG2} explains that the existence and continuity follow by descent from the scheme case; the adjunction follows similarly for smooth morphisms (cf. Section 4 of \cite{BNa}).  The identical arguments prove the assertions for twisted $\D$-modules.  
\end{proof}

\subsection{Microlocal Derived Categories}\label{sec:MDC}
Suppose that $M$ is a $\D$-module on $X$.  If $M$ is coherent, its singular support $SS(M)$ is a closed subset of $T^*X$.  For arbitrary (quasicoherent) $M$, we let $SS(M)$ denote the union of the closed subsets $SS(M')$ over all coherent $\D$-submodules $M'\subseteq M$. 
\subsubsection{Microlocal Quotient Category} 
For any open subset $U\subseteq T^*X$ with complement $S = T^*X\smallsetminus U$ we let $\D(X)_S$ denote the full DG subcategory of $\D(X)$ consisting of complexes whose cohomologies have singular support in $S$.  This subcategory is closed under taking colimits and direct summands.  By \cite{Keller, Drinfeld}, there is a quotient dg category
\bd
\D(U):= \D(X) / \D(X)_S.
\ed
\subsubsection{Microdifferential Operators}
Associated to the open set $U$ are two further categories.  

First, recall the sheaf $\cE_{X}$ of algebraic microdifferential operators on $T^*X$.  It is a sheaf of filtered algebras on $T^*X$, complete with respect to the inverse limit under the filtration; we write $\cE_X(0)$ for the subalgebra of elements of nonpositive degree.  Let $p_X: T^*X\rightarrow X$ denote the projection.  Then $\cE_X$ comes equipped with a faithfully flat homomorphism $p_X^{-1}\D_X \rightarrow \cE_X$.  For any conical open set $U\subseteq T^*X$ we write $\cE_U = \cE_X|_U$.  

There is a reasonable notion of coherent $\cE_U$-module and coherent $\cE_U(0)$-module.  A {\em lattice} for a  coherent $\cE_U$-module $M$ is a coherent $\cE_U(0)$-module $M(0)$ with
$\cE_U\cdot M(0) = M$.  A coherent $\cE_U$-module is {\em good} if it admits a lattice.    
\begin{notation}
We abusively write $\cE_U-\on{mod}$ for the abelian category of good $\cE_U$-modules.  We write $\cE_U-\on{Mod}$ for the ind-category of $\cE_U-\on{mod}$.
\end{notation}
The category $\cE_U-\on{Mod}$ is a (locally Noetherian) Grothendieck category.  

\subsubsection{$\cW$-Modules}
Again, let $U$ be a conical open set of $T^*X$.  The symplectic variety $T^*X$ comes equipped with a deformation quantization $\cW(0)$, a flat $\C[\![\hbar]\!]$-algebra that is $\Gm$-equivariant if $\Gm$ acts on $\hbar$ with weight $1$; see \cite{KR, BDMN} for further information.  Let $\cW = \cW(0)[\hbar\inv]$.  There is a natural homomorphism $\cE_X\rightarrow \cW$.    

The notions of coherent $\cW(0)$-module and coherent $\cW$-module, and their $\Gm$-equivariant analogues, are reviewed in \cite{BDMN}.  Recall that a {\em lattice} for a coherent $\Gm$-equivariant $\cW$-module $M$ is a coherent $\Gm$-stable $\cW(0)$-submodule $M(0)$ such that $M(0)[\hbar\inv] = M$.  The module $M$ is {\em good} if it admits a lattice.  For a conical open $U$, we write $\cW_U-\on{mod}$ for the abelian category of good 
$\Gm$-equivariant $\cW_U$-modules and $\cW_U-\on{Mod}$ for its ind-category.

\subsubsection{Equivalence of Different Types of Categories}
Let $\D-\on{mod}$ denote the abelian category of coherent $\D_X$-modules and $\D-\on{Mod}$ the abelian category of quasicoherent $\D_X$-modules.    
\begin{prop}
Suppose that $X$ is smooth and quasiprojective and $U\subseteq T^*X$ is a conical open set.  Let $S = T^*X\smallsetminus U$.  
Then:
\begin{enumerate}
\item The  natural functor $\D-\on{Mod}\rightarrow \cE_U-\on{Mod}$ factors through the quotient $\D-\on{Mod}/\D-\on{Mod}_S$.  
\item The resulting functors
\bd
\D-\on{Mod}/\D-\on{Mod}_S \longrightarrow \cE_U-\on{Mod} \longrightarrow \cW_U-\on{Mod}
\ed
are exact equivalences of abelian categories.  Moreover, they identify the images of $\D-\on{mod}$, $\cE_U-\on{mod}$, and $\cW_U-\on{mod}$.
\end{enumerate}
\end{prop}
\begin{proof}
See \cite{BDMN}.
\end{proof}
\begin{notation}\label{not:modcats}
We write $\cE(U)$ for the unbounded derived category of $\cE_U-\on{Mod}$ and $\cW(U)$ for the unbounded derived category of $\cW_U-\on{Mod}$.
\end{notation}
\begin{corollary}\label{cor:DorEorW}
The natural functors
$\D(U)\longrightarrow \cE(U)\longrightarrow \cW(U)$
are $t$-exact quasi-equivalences.
\end{corollary}

\subsection{Equivariant Microlocal Categories}\label{sec:EMC}
Suppose $\resol$ is a smooth symplectic variety with $\Gm$-action for which the symplectic form has weight $\ell$.  Suppose $\resol$ is equipped with a 
$\Gm$-equivariant $\cW$-algebra $\cW = \cW_{\resol}$.  As in the previous section, for any conical (i.e. $\Gm$-stable) open subset $U$ of $\resol$ we write
$\cW(U)$ for the unbounded derived category of the indization of the category of good $\Gm$-equivariant $\cW$-modules.  

Suppose $\gK$ is an affine algebraic group with a Hamiltonian action on $\resol$ and a fixed moment map $\mu^{\on{cl}}: \resol\rightarrow \mathfrak{k}^*$.  We assume that $\cW$ comes equipped with a {\em quantum comoment map} 
$\mu_c: \mathfrak{k}\rightarrow \cW(0)$
in the sense described in \cite{McN}.  On each product $(T^*\gK)^n\times \resol$, we have also a $\Gm^{n+1}$-equivariant $\cW$-algebra, $\cW_{\resol,n} := \cW_{(T^*\gK)^n\times \resol}$ given by inverting $\hbar$ in the completed tensor product  of $\cW_{\gK}(0)^{\otimes n}$ with $\cW_{\resol}(0)$.  

We now give the microlocal version of our simplicial description of the dg category of equivariant $\D$-modules. The key to this description is simply to note that the transfer bimodules which define the pull-back functors have a natural microlocal description and hence have analogues for $\cW$-algebras. (The general microlocal situation for $\mathcal E$-modules is dealt with in, for example, \cite[Chapter II, \S 3]{Sch}, here our description is somewhat more explicit using comoment maps.) The $\cW$-algebra $\cW_{T^*\gK}(0)$ comes with a canonical choice of quantum comoment map for $\mathfrak{k} = \on{Lie}(\gK)$.  For each $1\leq k\leq n$,  we let $\mu_k\colon \mathfrak{k}\to \cW_{\resol,n}(0)$ denote the composition of the quantum comoment map for $\cW_{T^*\gK}(0)$ with the natural inclusion of $\cW_{T^*\gK}(0)$ into $\cW_{\resol,n}(0)$ as the $k$th copy of $\cW_{T^*\gK}(0)$. We write $\mu_{0,1} = \mu_1$.  For each $k$ with $2\leq k \leq n$, let $\mu_{k-1,k} = \mu_{k-1}+\mu_k$ and define $B_{n,k}$ to be the quotient module $\cW_{\resol,n}/\cW_{\resol,n}\mu_{k-1,k}(\mathfrak{k})$. Similarly, let $B_{n,n}$ be the quotient module of $\cW_{\resol, n}$ by the ideal generated by the image of $(\mu_n+\mu_c)(\mathfrak{k})$. 

\begin{lemma}
Each quotient $B_{n,k}$ is naturally a $\big(\cW_{(T^*\gK)^n\times \resol}, a_{n,k}^{-1}\cW_{(T^*\gK)^{n-1}\times \resol}\big)$-bimodule, where
\bd
a_{n,k}: (T^*\gK)^n\times \resol \longrightarrow (T^*\gK)^{n-1}\times \resol
\ed
is the map induced on $(T^*\gK)^n\times \resol$ by the $k$th action map.  Moreover, $B_{n,k}$ is faithfully flat over $\cW_{(T^*\gK)^{n-1}\times \resol}$.
\end{lemma}
As in Notation \ref{not:modcats}, let $\cW((T^*\gK)^n\times \resol)$ denote the unbounded derived category of $\cW_{(T^*\gK)^n\times \resol}-\on{Mod}$. 
The following is then immediate:
\begin{prop}
The bimodules $B_{n,k}$ together with the projection define a cosimplicial object in dg categories:
\begin{equation}\label{eq:W-diagram}
\xymatrix{\cW(\resol) \ar@<0.5ex>[r]\ar@<-0.5ex>[r] &  \cW(T^*\gK\times \resol) \ar@<0.7ex>[r]\ar@<-0.7ex>[r]\ar[r]
& \cW(T^*\gK\times T^*\gK\times\resol) \ar@<0.9ex>[r] \ar@<0.3ex>[r] \ar@<-0.3ex>[r] \ar@<-0.9ex>[r] & \dots}.
\end{equation}
\end{prop}
\begin{defn}
The homotopy limit of the diagram \eqref{eq:W-diagram} is denoted by $\cW(\resol/\!\!/\gK,c)$.
\end{defn}

\subsection{Basic Properties}
Suppose next that $U\subseteq T^*X$ is a $G$-invariant open subset.  The following is immediate from Lemma \ref{holims of quotients}.
\begin{prop}\label{prop:holim}
Write $\D(U/\!\!/G,c) := \D(X/G,c)/\D(X/G,c)_S$.  Then:
\bd
\D(U/\!\!/G,c) \simeq \on{holim}\left(\xymatrix{\D(X)/\D(X)_S \ar@<0.5ex>[r]^{\hspace{-2em}a_c^!}\ar@<-0.5ex>[r]_{\hspace{-2em}p^!} &  \D(G\times X)/\D(G\times X)_{T^*G\times S} \ar@<0.7ex>[r]\ar@<-0.7ex>[r]\ar[r] & \dots }\right).
\ed
\end{prop}
It follows from Propositions \ref{prop:descent-of-adjoints} and \ref{prop:holim} and Corollary \ref{cor:DorEorW} that we can use $\cE$-modules or $\cW$-modules to study (or even define) the homotopy limit category $\D(U/\!\!/G,c)$.   
\begin{notation}
Suppose $U^\circ \subset U \subset T^*W$ are open complements to unions of KN strata.  We write
\begin{center}
$j: U^\circ \hookrightarrow U, \hspace{2em} \overline{j}: U^\circ \hookrightarrow T^*X$ \hspace{2em}
for the inclusions.
\end{center} 
\end{notation}
We abusively write $j^*$ to denote both the quotient functor $\D(U,c)\rightarrow \D(U^\circ,c)$ and the quotient functor $\D(U/\!\!/G,c)\rightarrow \D(U^\circ/\!\!/G,c)$.  These functors are automatically continuous.

\begin{prop}\label{jlowerstar}
The functors $j^*$ admit right adjoints $j_*$.  Moreover:
\mbox{}
\begin{enumerate}
\item Each adjoint $j_*$ is continuous and {\em cohomologically bounded}: there is some $N$ such that if $M\in \D(U^\circ/\!\!/G)^{[a,b]}$ then $j_*M\in\D(U/\!\!/G)^{[a,b+N]}$ and similarly for $\D(U^\circ)$.  
\item  $j^*\circ j_*\simeq \on{Id}$; in particular, $j_*$ is quasi-faithful.
\item The following diagram commutes:
\bd
\xymatrix{\D(U^\circ/\!\!/G)\ar[d]_{j_*} \ar[r]^{\pi^!}& \D(U^\circ) \ar[d]^{j_*}\\
\D(U/\!\!/G) \ar[r]^{\pi^!} & \D(U).
}
\ed
\end{enumerate}
\end{prop}
\begin{proof}
Existence of $j_*$ follows from \cite{Drinfeld} or \cite{Keller}.  The equality $j^*\circ j_*\simeq \on{Id}$ also follows from the construction.

To prove boundedness in (1),
it is enough to check after applying a faithful $t$-exact functor compatible with $\ast$-pushforward; hence it suffices to check after applying $\pi^!$.  Now apply \v{C}ech theory as in \cite[Theorem~7.14]{McNderived}.\ff{We note that the \v{C}ech theory of \cite{McNderived} uses microlocalizations.  However, in the case of $\D(X)$ (i.e., working on the flat cover rather than the equivariant category $\D(X/G,c)$), one can use localizations of the rings instead and arrive at an equivalent \v{C}ech theory that works more naively for quasi-coherent $\D_X$-modules.} This also yields continuity.

Commutativity of the diagram in (3) follows from Proposition \ref{prop:descent-of-adjoints}.
\end{proof}
The category $\D(X/G,c)_{\on{coh}}$ consists of those complexes with only finitely many nonzero cohomologies, each of which is a coherent $(G,c)$-equivariant $\D$-module.  
Its image in $\D(U^\circ/\!\!/G,c)$ under $j^*$ where $j:U^\circ\hookrightarrow T^*X$ is the inclusion is, by definition, $\D(U^\circ/\!\!/G,c)_{\on{coh}}$.  Similarly, $\D(X/G,c)^c$ denotes the subcategory of compact objects of $\D(X/G,c)$ and $\D(U^\circ/\!\!/G,c)^c$ denotes the subcategory of compact objects of $\D(U^\circ/\!\!/G,c)$.  

\begin{lemma}\label{lem:compact is coh}
One has $\D(X/G,c)^c\subseteq \D(X/G,c)_{\on{coh}}$ and $\D(U^\circ/\!\!/G,c)^c\subseteq \D(U^\circ/\!\!/G,c)_{\on{coh}}$.
\end{lemma}
\begin{proof}
We use that $\pi^!$ has a continuous right adjoint $\pi_*$ (Lemma \ref{lem:pi-functors}) and that compactness implies coherence in $\D(X)$ and $\D(U^\circ)$.  The conclusion is immediate from Proposition \ref{lem:cptimpliescoh}.
\end{proof}

\subsection{Verdier Duality}
Suppose $X$ is a smooth scheme.  The Verdier duality functor $\vD: \D(X)_{\on{coh}}^{op}\rightarrow \D(X)_{\on{coh}}$ defined by $\vD(M) = \Hom_{\D}(M, \D\otimes\omega)$ is an anti-self-equivalence of the coherent derived category.  Applying Verdier duality level-wise for a simplicial scheme defines a Verdier duality anti-equivalence of the subcategory of coherent objects on the corresponding algebraic stack as in Section 4.1 of \cite{BNa}.  In the case of the twisted diagram \eqref{twisted-diagram}, Verdier duality does {\em not} preserve the twist $c$: it takes the category for twist $c$ to the category for another twist $c^{\vee}$ (where $c^{\vee} = \rho-c$; however, this will not be important for us). 
\begin{prop}
Let $S$ be a $G$-invariant conical closed subset of $T^*X$.  Then the full subcategory $\big(\D(X/G, c)_S\big)_{\on{coh}}$ of $\D(X/G,c)_{\on{coh}}$ is mapped by Verdier duality to the corresponding subcategory $\big(\D(X/G, c^{\vee})_{\on{coh}}\big)_S$ of $\D(X/G,c^{\vee})_{\on{coh}}$.
\end{prop}
The statement of the proposition is standard for $\D(X)$: see \cite[Proposition~D.4.2]{HTT}.  It follows for $\vD$ on $\D(X/G,c)$ since $\vD$ is defined there by descent.
\begin{corollary}
There exists a Verdier duality anti-equivalence $\vD: \D(U/\!\!/G,c)_{\on{coh}}^{\on{op}}\rightarrow \D(U/\!\!/G, c^{\vee})_{\on{coh}}$ such that $\vD\circ j^* \simeq j^*\circ \vD$ on $\D(X/G,c)_{\on{coh}}$.
\end{corollary}
A crucial property of the functor $\vD$ is that it preserves compactness: 
\begin{lemma}\label{lem:dualcompact}
\mbox{}
\begin{enumerate}
\item  $\vD\big(\D(X/G,c)^c\big)=\D(X/G,c^\vee)^c$. 
\item $\vD\big(\D(U/\!\!/G,c)^c\big) = \D(U/\!\!/G,c^\vee)^c$.
\end{enumerate}
\end{lemma}
\begin{proof}
The first statement is Corollary 8.4.2 of \cite{DG2}, at least in the untwisted case; the twisted case follows exactly the same argument.  For (2), observe that 
the functor $j_*$ is faithful and, by Proposition \ref{jlowerstar}, continuous; hence its left adjoint preserves compact objects (cf. \cite{Neeman} or \cite[Proposition~1.2.4]{DG}). Moreover, the subcategory $\D(U/\!\!/G,c)^c$ of compact objects is Karoubi-generated by the images of compact objects of $\D(X/G,c)$ by Lemma \ref{lem:karoubigenerate}.  Hence Verdier duality preserves compactness in $\D(U/\!\!/G,c)$ by construction of $\vD$.  
\end{proof}

 \subsection{A Version of the Projection Formula}
 Recall that the projection formula for $\D$-modules for a morphism of smooth schemes $f: X\rightarrow Y$ says that, if $M\in\D(X)$ and $N\in\D(Y)$, then, in the notation of \cite{DG2},
\bd
f_{\on{dR},*}(f^!(N)\otimes M) \simeq N\otimes f_{\on{dR},*}M.
\ed
  Applying the functor $\Gamma_{\on{dR}}$ of global de Rham cohomology on $Y$ gives an isomorphism that can be written more algebraically as 
  \bd
  \Gamma(X, f^!N \otimes_{\D_X} M) \simeq \Gamma(Y, N\otimes_{\D_Y} f_{\on{dR},*} M);
  \ed
   here the two tensor products are naturally only (complexes of) sheaves of vector spaces, and $\Gamma(X,-)$ and $\Gamma(Y,-)$ denote usual sheaf (hyper)cohomology.   The same formulas hold for twisted $\D$-modules.
 
 \begin{comment}
 As in \cite{DG2}, giving conditions under which similar formulas hold for a morphism of (derived) stacks is more complicated; the notion of {\em safety} plays a role, cf. Section 8.2 of \cite{DG2}.  We will prove the following version of a projection formula, which will suffice for our purposes.  As before, assume that $\openset$ is open in $\bigvar=T^*(W/G)$ and keep the notation of earlier sections.  
 \end{comment}
 
\begin{prop}\label{projection formula}
Let $j: U^\circ\hookrightarrow T^*X$ be the open complement to a union of KN strata.  
  Then for any $M\in \D(X/G,c^\vee)$ and $N\in \D(X/G,c)$, $j_*j^*M\otimes_{\D} N \simeq M\otimes_{\D} j_*j^*N$.  
\end{prop}
 \begin{proof}
Let $\wt{M} = \pi^!M$ and $\wt{N} = \pi^!N$.  We get commutative diagrams
\bd
\xymatrix{j_*j^*\wt{M}\otimes_{\D} \wt{N} \ar[d]_{\simeq} \ar[d]_{\simeq} \ar@<0.5ex>[r]^{a_c^!}\ar@<-0.5ex>[r]_{p^!} & p^!j_*j^*\wt{M}\otimes p^!\wt{N}\ar[d]^{\simeq} \ar@<0.7ex>[r]\ar@<-0.7ex>[r]\ar[r] & \dots\\
\wt{M}\otimes_{\D} j_*j^*\wt{N} \ar@<0.5ex>[r]^{a_c^!}\ar@<-0.5ex>[r]_{p^!} & p^!\wt{M}\otimes_{\D} p^!j_*j^*\wt{N} \ar@<0.7ex>[r]\ar@<-0.7ex>[r]\ar[r] & \dots,
}
\ed
where the vertical arrows are isomorphisms from the \v{C}ech complex argument of \cite[Theorem~6.3]{McNderived} (it is essentially Formula (6.1)). 
Since $-\otimes_{\D} -$ is defined on $X/G$ via the homotopy limit \eqref{twisted-diagram}, the statement of the proposition follows.
\end{proof}
\begin{corollary}\label{coherence is enough-concrete}
Let $j: U^\circ\hookrightarrow T^*X$ be the open complement to a union of KN strata. 
Suppose that $j_*$ preserves coherence for all $c$: that is, 
\bd
\text{if}\,\,\, M \in \D(U^\circ/\!\!/G,c)_{\on{coh}} \,\, \text{then}\,\, j_*M \in \D(X/G,c)_{\on{coh}}.
\ed
Then $j_*$ preserves compactness for all $c$ as well.  
\end{corollary}
\begin{proof}
Apply Proposition \ref{coherence is enough}, using Proposition \ref{projection formula} to obtain hypothesis (3).
\end{proof}

 \section{Proof of Theorem \ref{main theorem}}
 We begin with the following slight revision of Theorem \ref{main theorem} (the remaining properties of Theorem \ref{main theorem} are immediate).
 \begin{thm}\label{main theorem-alt}
\mbox{}
\begin{enumerate}
\item
The functor $j^*$ admits a left adjoint $j_!$.  
\item The left adjoint $j_!$ takes compact objects to compact objects.
\item The right adjoint $j_*$ of $j^*$ takes bounded coherent complexes to bounded coherent complexes.
\item The functor $i^!$ admits a continuous right adjoint.
\item The functor $i_*$ admits a left adjoint $i^*$.
\item The left adjoint $i^*$ takes compact objects to compact objects.  
\item The right adjoint $i^!$ takes bounded coherent complexes to bounded coherent complexes.

\end{enumerate}
\end{thm}
 \subsection{Strategy of Proof}
 We begin by explaining the strategy of the proof of Theorem \ref{main theorem-alt}.  
We note that all the statements of the two Theorems reduce to Zariski-local statements; hence we assume that $\mathcal{Y}=X/G$ is a quotient of a smooth variety by a reductive group.
 
Consider the inclusions of open sets complementary to unions of KN strata,
\bd
j: U^\circ \hookrightarrow U, \hspace{1em} \overline{j}: U^\circ \hookrightarrow T^*X.
\ed
To prove the theorem, by Proposition \ref{basic equivalences}, it suffices to show that $j_*$ preserves compactness.  We note that, to prove that 
$\overline{j}_*$ preserves compactness, it suffices, by Proposition \ref{coherence is enough-concrete}, to prove:
\begin{prop}\label{prop:fullcoherence}
The functor $\overline{j}_*$ preserves coherence. 
\end{prop}
\noindent
The proof of Proposition \ref{prop:fullcoherence} appears in Section \ref{proofoffullcoherence}.  By an obvious induction, it suffices to consider the case when $U\smallsetminus U^\circ$ consists of a single KN stratum $S_\beta$; moreover, choosing an appropriate refinement of the partial order on KN strata to a total order, we may assume that 
\bd
U = T^*X\smallsetminus S_{>\beta} \hspace{2em} \text{and} \hspace{2em} U^\circ = T^*X\smallsetminus S_{\geq\beta}.
\ed
  We assume this for the remainder of the proof.

Propositions \ref{prop:fullcoherence} and \ref{prop:presofcompactness} together establish assertions (1) through (4) of Theorem \ref{main theorem-alt}.  The remainder are established in Section \ref{sec:therest}.
 
\subsection{Pullback to Slices}
Recall that we have fixed a KN stratum $S_\beta$.  Let $\gK = U^-_\beta\ltimes \Gm$ as in Section \ref{sec:sympl-geom}.  For each $n$, there is a natural pullback functor on derived categories of $\D$-modules $\D(G^n\times X)\rightarrow \D(\gK^n\times X)$.  These functors assemble into a functor of cosimplicial diagrams
\bd
\xymatrix{\D(X)\ar[d] \ar@<0.5ex>[r]^{\hspace{-1em}a_c^!}\ar@<-0.5ex>[r]_{\hspace{-1em}p^!} &  \D(G\times X)\ar[d] \ar@<0.7ex>[r]\ar@<-0.7ex>[r]\ar[r]
& \D(G\times G \times X) \ar[d]\ar@<0.9ex>[r] \ar@<0.3ex>[r] \ar@<-0.3ex>[r] \ar@<-0.9ex>[r] & \dots \\
\D(X) \ar@<0.5ex>[r]^{\hspace{-1em}a_c^!}\ar@<-0.5ex>[r]_{\hspace{-1em}p^!} &  \D(\gK\times X) \ar@<0.7ex>[r]\ar@<-0.7ex>[r]\ar[r]
& \D(\gK\times \gK \times X) \ar@<0.9ex>[r] \ar@<0.3ex>[r] \ar@<-0.3ex>[r] \ar@<-0.9ex>[r] & \dots}
\ed
 yielding a functor $\D(X/G,c)\rightarrow \D(X/\gK,c)$.  This functor is well known to be faithful and $t$-exact.
 
We now consider the \'etale map 
 \begin{equation}\label{eq:cover}
 \ds \coprod_{z\in Z_\beta^{ss}} U^-\times U_{D_z} \rightarrow \left(T^*X\smallsetminus S_{>\beta}\right)
 \end{equation}
 from Section \ref{sec:sympl-geom}.  This map is $\gK$-equivariant.  Since it is \'etale, the Kontsevich quantizations associated to the Poisson structures $\hbar\{\bullet,\bullet\}$ are compatible with the pullback.  For each $n$, we get a pullback functor 
 $a^*: \cW\big((T^*\gK)^n\times T^*X\big)\rightarrow \cW\big(T^*\gK)^n\times \coprod_{z\in Z_\beta^{ss}} U^-\times U_{D_z}\big)$.  These functors are compatible with the corresponding cosimplicial diagrams of dg categories and induce a pullback functor
$a^*: \cW(T^*X/\!\!/\gK, c)\rightarrow \prod_z \cW\big(U^-\times U_{D_z}/\!\!/\gK, c\big)$, where we note that the product on the right is finite since we have chosen finitely many $z$ such that the images of the corresponding maps cover an open set in $T^*X\smallsetminus S_{>\beta}$ whose $G$-orbit covers $S_\beta$.
 This latter functor generally need not be faithful.  
 
\begin{notation}
We let 
 $V\subset T^*X$ denote the $G$-orbit of the image of \eqref{eq:cover} and let $V^\circ = V\cap U^\circ$.
 \end{notation}
  Both $V^\circ$ and $V$ are $G$-stable open sets.
 We have:
 \begin{prop}\label{prop:faithful-t-exact}
 The composite functor 
 \bd
\xymatrix{
 \D(V/\!\!/G,c)\ar[r] \ar@/^2pc/[rr]^{\mathsf{C}^*}&  \D(V/\!\!/\gK,c)\ar[r]^{\hspace{-2em} a^*} & \prod_{z} \cW\big(U^-\times U_{D_z}/\!\!/\gK,c\big)
 }
 \ed
 is faithful and $t$-exact.
 \end{prop}
 \begin{proof}
 The corresponding maps of schemes are all flat, so it follows that the functor is $t$-exact.  Since the diagram
 \begin{equation}\label{eq:pullback-square}
 \begin{gathered}
 \xymatrix{\D(V/\!\!/G,c) \ar[d]\ar[r]^{\hspace{-2em} \mathsf{C}^*} & \prod_{z} \cW\big(U^-\times U_{D_z}/\!\!/\gK,c\big)\ar[d]\\
 \D(V) \ar[r]^{\hspace{-2em} a^*} & \prod_{z} \cW\big(U^-\times U_{D_z}\big)}
 \end{gathered}
 \end{equation}
  commutes, it suffices to check that the composite $\D(V/\!\!/G,c) \rightarrow \D(V) \rightarrow \prod_{z} \cW\big(U^-\times U_{D_z}\big)$ is faithful.  This follows for $\Gm$-equivariant $\cW(0)$-modules, by passing to $\hbar=0$, from Corollary \ref{cor:fflat}; hence it follows for $\Gm$-equivariant $\cW$-modules.
 \end{proof}

 \subsection{Proof of Proposition \ref{prop:fullcoherence}}\label{proofoffullcoherence}
 We prove Proposition \ref{prop:fullcoherence} by proving that the $*$-pushforward across a single closed KN stratum preserves coherence.  Since the KN strata are partially ordered by closure relations, the full statement immediately follows by an induction.  
\subsubsection{Step 1: Reduction to Preservation of Coherence on $U^-\times U_{D_z}$.}
For each $z\in Z_\beta^{ss}$ and choice of $U_{D_z}$, we get a closed subset $\wt{Y}_\beta\subset U^-\times U_{D_z}$, the set of points that converge under the downward $\beta$-flow to the $\beta$-fixed locus.  We have $a(\wt{Y}_\beta)\subset S_\beta\subseteq T^*X\smallsetminus S_{>\beta}$, and in fact:
\begin{lemma}
We have $\wt{Y}_\beta \cong Y_\beta$ and $a^{-1}(S_\beta) = U^-\times \wt{Y}_\beta$. 
\end{lemma}
\begin{proof}
Since $\beta(\Gm)$ has negative weights on $U^-$, we get $\wt{Y}_\beta   = a\inv(Y_\beta) \cong Y_\beta$.  Thus $a\inv(U^-\cdot Y_\beta) = U^-\times \wt{Y}_\beta$.  Now by construction, $a$ is quasi-finite \'etale, so 
$a^{-1}(S_\beta)$ is smooth and equidimensional; hence $a\inv(U^-\cdot Y_\beta)$ is dense in $a\inv(S_\beta)$.
 It follows that $a\inv(S_\beta\smallsetminus Y_\beta)$ is a subset of $U^-\times U_D$ contained in the closure of $U^-\times Y_\beta$.  But the latter is already closed.
\end{proof}
 Let $\wt{j}: \big(U^-\times U_{D_z}\big)\smallsetminus \wt{Y}_\beta \hookrightarrow U^-\times U_{D_z}$ denote the open immersion, and 
\bd
\wt{j}_*: \cW\big((U^-\times U_{D_z})\smallsetminus \wt{Y}_\beta/\!\!/\gK, c\big) \rightarrow \cW(U^-\times U_{D_z}/\!\!/\gK, c)
\ed
the corresponding adjoint to $j^*$.  
 \begin{lemma}\label{lem:base-change}
 We have $\wt{j}_*\mathsf{C}^* = \mathsf{C}^*j_*$.
 \end{lemma}
 \begin{proof} 
 Both functors $\wt{j}_*$ and $j_*$ are defined via Proposition \ref{prop:descent-of-adjoints}, so it suffices to check the claim level-wise on the cosimplicial diagrams of categories used to define the two equivariant categories.  Here the statement follows since both $\wt{j}_*$ and $j_*$ are calculated by \v{C}ech complexes, which are compatible both with flat base change for $\D$-modules and with \'etale base change for $\cW$-modules.
 \end{proof}
 
 \begin{prop}\label{prop:descentofcohpres}
 Suppose that, for all $z\in Z_\beta^{ss}$ and each \'etale chart $U^-\times U_D\rightarrow T^*X$ around $z$, the functor $\wt{j}_*$ preserves coherence.
 Then $j_*$ preserves coherence.
 \end{prop}
 \begin{proof}
 First, note that preservation of coherence for $j_*$ can be checked Zariski-locally.  Hence it suffices to check for the functor $j_*$ corresponding to the inclusion 
 $j: V^\circ\rightarrow V$ coming from the covering \eqref{eq:cover}.

 Taking a finite product as in \eqref{eq:pullback-square}, consider the functor $\mathsf{C}^*$; it is faithful and $t$-exact by Proposition
 \ref{prop:faithful-t-exact}.  Hence if $\mathsf{C}^*(M)$ is coherent, so is $M$.  Suppose $M\in \D(V^\circ/\!\!/G,c)$ is coherent.  Then $\mathsf{C}^*M$ is coherent.  By Lemma \ref{lem:base-change} and the assumption, $\wt{j}_*\mathsf{C}^*M = \mathsf{C}^*j_*M$ is coherent, hence so is $j_*M$.   
 \end{proof}
 \subsubsection{Step 2: Proof of Preservation of Coherence for $\wt{j}_*$.}
 We have thus reduced to proving that, for a single $U^-\times U_{D_z}$, letting $\wt{j}: U^-\times U_{D_z} \smallsetminus \wt{Y}_\beta \hookrightarrow U^-\times U_{D_z}$, that 
 \bd
 \wt{j}_*: \cW\big(((U^-\times U_{D_z}) \smallsetminus \wt{Y}_\beta)/\!\!/\gK, c\big) \rightarrow \cW\big((U^-\times U_{D_z})/\!\!/\gK, c\big)
 \ed
 preserves coherence.  By Proposition \ref{basic equivalences}, to do this, it suffices to prove that $\wt{j}_!$ is defined on a collection of compact generators.  
 
 From now on we write $U_D := U_{D_z}$.  
 Recall the \'etale morphism $\phi: U^-\times U_D\rightarrow T^*U^-\times \cS$ from \eqref{eq:proj}.  As in \cite{McN}, this map is compatible with symplectic structures and $\gK$-equivariant.  Moreover, standard arguments show that 
$\cW\big((T^*U^-\times \cS)/\!\!/\gK, c\big)\simeq \cW(\cS/\!\!/\beta(\Gm), c)$.  We thus get a pullback functor $\phi^*:\cW(\cS/\!\!/\beta(\Gm), c)\rightarrow 
\cW\big((U^-\times U_D)/\!\!/\gK, c\big)$.  Moreover, since $\phi(\wt{Y}_\beta)$ equals an open subset of the corresponding locus $\overline{Y}_\beta$ (of points that converge under the downward $\beta$-flow to the $\beta$-fixed locus) in $\cS$ and in fact $\wt{Y}_\beta = \phi\inv(\overline{Y}_\beta)$ (see \cite{McN}), we get also a 
functor 
\bd
(\phi^\circ)^*:\cW\big((\cS\smallsetminus \overline{Y}_\beta)/\!\!/\beta(\Gm), c\big)\longrightarrow 
\cW\big(((U^-\times U_D)\smallsetminus \wt{Y}_\beta)/\!\!/\gK, c\big).
\ed

We fix the identity  character $\chi: \Gm\rightarrow \Gm$, $\chi(z) = z$ (this choice is for notational convenience when referring to \cite{McNderived, McN}).   We then define objects $M_c(\chi^\ell)$ of $\cW\big(\cS/\!\!/\beta(\Gm), c\big)^\heartsuit$, quantizing line bundles, following Formula (4.2) of Section 4.1  of \cite{McNderived}.  More specifically, fixing the quantum comoment map $\mu$ associated to the infinitesimal $\mathfrak{g}$-action, for a character 
$a: \mathfrak{g}\rightarrow\C$ we write $\mu_a = \mu+a$.   Write $\cS = T^*Q$, so that the deformation quantization $\theo^{\hbar}(\cS)$ is naturally identified with a completion of $\D(Q)$.  Let  $M_{c-\ell d\chi} = \D/\D\mu_{c-\ell d\chi}\big(\on{Lie}(\Gm)\big)$; as in \cite{McNderived}, the latter is a $(c-\ell d\chi)$-twisted $\Gm$-equivariant $\D$-module on $Q$, hence naturally defines an object of $\cW\big(\cS/\!\!/\beta(\Gm), c-\ell d\chi\big)^\heartsuit$.  
Then, as in Formula (4.2) of \cite{McNderived}, 
$M_c(\chi^\ell) = M_{c-\ell d\chi} \otimes \chi^\ell$.  
\begin{prop}\label{prop:cpt-gen-on-S}
\mbox{}
\begin{enumerate}
\item The objects $M_c(\chi^\ell)[N]$ are compact.
\item
For any choice of $\ell_0$, the restrictions of the objects $M_c(\chi^\ell)[N]$ ($\ell\in\mathbb{Z}_{\leq \ell_0}, N\in\mathbb{Z}$) to $\cS\smallsetminus \overline{Y}_\beta$ constitute a collection of compact generators of $\cW\big((\cS\smallsetminus \overline{Y}_\beta)/\!\!/\beta(\Gm), c\big)$.
\item There exists a $d\geq 1$ such that for any $k$, the full subcategory of $\cW\big((\cS\smallsetminus \overline{Y}_\beta)/\!\!/\beta(\Gm), c\big)$ K-generated
by the restrictions of $M_c(\chi^k), M_c(\chi^{k+1}), \dots, M_c(\chi^{k+d})$ contains $M_c(\chi^\ell)$ for all $\ell\leq k$.  
\end{enumerate}
\end{prop}
\begin{proof}
For (1) and (2), repeat the argument of Proposition 4.14 of \cite{McNderived}. 

Assertion (3) is a direct analogue of generation statements like Theorem 4 of \cite{Orlov} (or, in a quantized setting, \cite{MV, Stadnik}).   Namely, note that 
$\overline{Y}_\beta\subset\cS$ is a linear subspace cut out by $\beta$-semi-invariants in $\C[\cS]$ of positive degree.  Choose a minimal list $f_1,\dots, f_m$ of such semi-invariants, with weights $d_1,\dots, d_m$, and let $L$ be the graded vector space they span.  Let $d= d_1 + \dots + d_m-1$.   This will be the choice of $d$ in assertion (3).   
\begin{remark}\label{remark:characters}
As a $\Gm$-representation, each vector space $\bigwedge^k L^*$ is a direct sum of characters $\chi^{-d_S}$ where $S\subset \{1,\dots, m\}$ is a $k$-element subset and $d_S = \sum_{i\in S} d_i$.
\end{remark}

Choose an identification $\cS = T^*Q$ (as in the description of $M_c(\chi^\ell)$ above) and lift $L$ to a semi-invariant subspace $L\subset \D(Q)$; note that the elements of $L$ pairwise commute as elements of $\D(Q)$ since $\beta$ preserves the symplectic structure on $\cS$ and the semi-invariants in $L$ all have positive weight.  We obtain a Koszul-type complex 
\bd
0 \rightarrow \D(Q) \otimes \bigwedge^m L^* 
\longrightarrow
\D(Q) \otimes \bigwedge^{m-1} L^* \longrightarrow \dots
\D(Q)\otimes \bigwedge L^* \longrightarrow \D(Q),
\ed
which is exact except at the right-hand end.  Twisting this complex by $\chi^N$ for some choice of $N$ and applying the ``$\Gm$-equivariantization functor'' $\Phi_c$ of Section 4.7 of \cite{McNderived}, we obtain a complex which is exact except at the right-hand end and whose terms are of the form $M_c(\chi^{d_S+N})$, where $d_S$ is as in Remark \ref{remark:characters}.  Restricting to a complex in $\cW\big((\cS\smallsetminus \overline{Y}_\beta)/\!\!/\beta(\Gm), c\big)$, we obtain an exact complex (since the right-hand cohomology has singular support in $\overline{Y}_\beta$ by construction) whose left-hand term is 
$M_c(\chi^{-d-1+N})$, whose right-hand term is $M_c(\chi^N)$, and whose intermediate terms are direct sums of objects $M_c(\chi^{N'})$ for 
$-d-1+N<N'<N$.   Taking $N = k+d$, it follows that any full Karoubi-closed subcategory that contains $M_c(\chi^k), \dots, M_c(\chi^{k+d})$ also contains 
$M_c(\chi^{k-1})$.  A descending induction now proves assertion (3). 
\end{proof}
 
 Next, using the induction functor of $\cW$-modules that quantizes \eqref{eq:proj} as in Section 6.4 of \cite{McN}, we obtain corresponding objects of
 $\cW\big(U^-\times U_D)/\!\!/\gK, c\big)$.  We denote the object of 
 $\cW\big(U^-\times U_D)/\!\!/\gK, c\big)$
 corresponding to $M_c(\chi^\ell)$ by $\theo_{(U^-\times U_D)/\!\!/\gK}^\hbar(\chi^\ell)$.  
 \begin{prop}\label{prop:cpctness-and-vanishing}
 \mbox{}
 \begin{enumerate}
 \item
 Each object $\theo_{(U^-\times U_D)/\!\!/\gK}^\hbar(\chi^\ell)$ is compact.  
 \item The restrictions of the objects 
 $\theo_{(U^-\times U_D)/\!\!/\gK}^\hbar(\chi^\ell)[N]$, ($\ell\in\mathbb{Z}_{\leq \ell_0}, N\in\mathbb{Z}$) to $U^-\ltimes U_D\smallsetminus \wt{Y}_\beta$  constitute a collection of compact generators of
  $\cW\big(((U^-\times U_D)\smallsetminus \wt{Y}_\beta)/\!\!/\gK, c\big)$.  
  \item There exists $k$ such that for $d$ as in Proposition \ref{prop:cpt-gen-on-S}(3), $k\leq \ell\leq k+d$, and all objects $M$ of $\cW\big(U^-\times U_D)/\!\!/\gK, c\big)$ with cohomologies supported in $\wt{Y}_\beta$, we have
  \bd
  \Hom(\theo_{(U^-\times U_D)/\!\!/\gK}^\hbar(\chi^\ell), M) \simeq 0.
  \ed
  \end{enumerate}
  \end{prop}
  \begin{proof}
  Assertions  (1) and (2) follow from Proposition \ref{prop:cpt-gen-on-S} since the induction functor has a faithful right adjoint (its faithfulness follows because the map \eqref{eq:proj} is \'etale).  Assertion (3) follows for objects $M$ of $\cW\big(U^-\times U_D)/\!\!/\gK, c\big)^\heartsuit$ by the vanishing Theorem 6.6 of \cite{McN}: indeed, we get
  \bd
  \begin{split}
  \Hom_{\cW\big(U^-\times U_D)/\!\!/\gK, c\big)}(\theo_{(U^-\times U_D)/\!\!/\gK}^\hbar(\chi^\ell), M) = 
  \Hom_{\cW\big(U^-\times U_D)/\!\!/\gK, c-\ell d\chi\big)}(\theo_{(U^-\times U_D)/\!\!/\gK}^\hbar\otimes \chi^\ell, M) = \\
  \Hom_{\cW\big(U^-\times U_D)/\!\!/\gK, c-\ell d\chi\big)}(\theo_{(U^-\times U_D)/\!\!/\gK}^\hbar, M\otimes \chi^{-\ell}).
  \end{split}
  \ed
 Theorem 6.6 of \cite{McN} implies that, for each $c$, there is an infinite consecutive subset of $\ell\in \mathbb{Z}$ for which the last Hom above vanishes for all objects $M$ of $\cW\big(U^-\times U_D)/\!\!/\gK, c-\ell d\chi\big)^\heartsuit$ supported in $\wt{Y}_\beta$.  
  A spectral sequence argument then implies the general statement of (3) (i.e. for arbitrary complexes with cohomologies supported in $\wt{Y}_\beta$).  
  \end{proof}
  \begin{corollary}
  \mbox{}
  \begin{enumerate}
  \item The functor $j_!$ is defined for a collection of compact generators of $\cW\big(U^-\times U_D)/\!\!/\gK, c\big)$.
  \item 
  The functor $\wt{j}_*$ preserves coherence.
  \end{enumerate}
  \end{corollary}
  \begin{proof}
  It follows from Proposition \ref{prop:cpctness-and-vanishing}(3) and 
  Proposition \ref{prop:shriekdefined} that $j_!$ is defined for the objects $\theo_{(U^-\times U_D)/\!\!/\gK}^\hbar(\chi^\ell)$, $k\leq \ell\leq k+d$.  It thus follows that $j_!$ is defined for the objects $\theo_{(U^-\times U_D)/\!\!/\gK}^\hbar(\chi^\ell)$ where $\ell\leq k+d$ by 
Proposition  \ref{prop:cpt-gen-on-S}(3) and the construction of those objects.  These form a collection of compact generators by 
Proposition \ref{prop:cpctness-and-vanishing}(2), thus proving assertion (1).
  Assertion (2) then follows by Proposition \ref{basic equivalences}.
  \end{proof}
 
As explained above, by induction this proves preservation of coherence for $j_*$.  Assertions (1), (2), (3), and (4) of Theorem \ref{main theorem-alt} follow by applying Lemma \ref{lem:jbarimpliesj} and Propositions \ref{prop:shriekdefined} and \ref{basic equivalences}. \hfill\qedsymbol

 \subsection{Proof of Theorem \ref{main theorem-alt}}\label{sec:therest}
 It remains to establish parts (5), (6), and (7) of Theorem \ref{main theorem-alt}.  Assertion (5) follows from (1) by Proposition \ref{prop:iupperstar}.  Assertion (6) follows from assertion (2) using the exact triangle of Proposition \ref{prop:iupperstar}.  Assertion (7) follows from assertion (3) using 
 Proposition \ref{prop:quotient-properties}(4).

 \section{Hochschild Homology and the Kirwan Map}\label{sec:Kirwan}
 In this section we briefly explain how Theorem \ref{thm:Kirwan} follows from Theorem \ref{main theorem}.
\subsection{Hochschild Homology of DG Categories}
 We begin with a categorically-oriented summary of Hochschild homology; a good reference is \cite{BNtraces}.  
 Consider a  dualizable DG category over $\C$.  Associated to each such category $\cC$ is a DG vector space $HH_*(\cC)$, the {\em Hochschild homology}, or in the terminology of \cite{BNtraces} the {\em trace} $\on{Tr}(\cC)$, of $\cC$.  Dualizability of the category means \cite[Section~2]{GaitsgoryDG} that there are another DG category $\cC^\vee$ and unit and co-unit morphisms $\on{dgVect}\rightarrow \cC^\vee\otimes\cC$, 
 $\cC^\vee\otimes\cC\rightarrow \on{dgVect}$ satisfying standard identities.  The Hochschild homology is then a dg vector space representing the 
 composite
$ \on{dgVect}\longrightarrow \cC^\vee\otimes\cC \longrightarrow \on{dgVect}.$
 
 Suppose next that $\cC\longrightarrow \cD$ is a continuous functor, i.e., a functor with a right adjoint.  Then, as explained in \cite[Proposition~1.12]{BNtraces}, there is an induced morphism on Hochschild homology, $HH_*(\cC)\xrightarrow{HH(f)} HH_*(\cD)$; such morphisms, moreover, are compatible with compositions of continuous functors.  
 
 In particular, the following is immediate from Theorem \ref{main theorem}:
 \begin{prop}
 Let $\cY$ be a smooth, quasicompact stack that is exhausted by quotient stacks.  Equip $T^*\cY$ with a Kirwan-Ness stratification and use notation as in Section \ref{sec:main-results}.  Then the homomorphism 
 \bd
 HH_*\big(\cD(\cU)\big)\xrightarrow{HH(j^*)} HH_*\big(\D(\cU^\circ)\big)
 \ed
 is surjective on cohomology groups.  
 \end{prop}
 \begin{proof}
 The functors $j_!$ and $j^*$ are continuous, and $j^*\circ j_!\simeq \on{Id}_{\D(\cU^\circ)}$.   We thus get homomorphisms
 \bd
HH_*\big(\D(\cU^\circ)\big) \xrightarrow{HH(j_!)}  HH_*\big(\cD(\cU)\big)\xrightarrow{HH(j^*)} HH_*\big(\D(\cU^\circ)\big)
  \ed
  whose composite is the identity.
 \end{proof}
 
\subsection{Proof of Theorem \ref{thm:Kirwan}}
The proof of Theorem \ref{thm:Kirwan} now rests on two known facts.  First:
\begin{prop}[\cite{BNtraces}]\label{HH of I}
Suppose $X$ is a smooth quasicompact variety and $G$ is an algebraic group acting on $X$.  Then $HH_*\big(\D(X/G)\big) \simeq H^*_{dR}\big(I(X)/G\big)$.
\end{prop}
The second fact requires a bit of background.  Fix notation as in Theorem \ref{thm:Kirwan}, assuming, in particular, that $G$ acts on $\mu\inv(0)^{ss}$ with free orbits, yielding a smooth symplectic quotient variety $\mu\inv(0)^{ss}/G$.  
Note that the $G$-equivariant (and $\Gm$-equivariant, for the scaling action) $\cW$-algebra $\cW_{T^*X}$ comes equipped with a quantum comoment map and descends as in \cite[Section~2.5]{KR} to a $\Gm$-equivariant $\cW$-algebra $\cW_{\mu\inv(0)^{ss}/G}$ on $\mu\inv(0)^{ss}/G$.  
Following Sections \ref{sec:MDC} and \ref{sec:EMC} and \cite[Proposition~2.8]{KR}, the category $\D(T^*X^{ss}/G)$ is equivalent to the (unbounded) derived category 
$D(\cW_{\mu\inv(0)^{ss}/G}-\on{Mod})$ (notation as in Section \ref{sec:MDC}).  

We now consider the category $\cW_{\mu\inv(0)^{ss}/G}-\on{Mod}$, the ind-category of the category of good $\Gm$-equivariant $\cW_{\mu\inv(0)^{ss}/G}$-modules.  Its Hochschild homology is the same as the Hochschild homology of its derived category.  

The following is a minor variant of the results of 
Section 2 of \cite{BrGe}.
\begin{prop}\label{HH of DQ}
We have
$HH_*\big(\cW_{\mu\inv(0)^{ss}/G}-\on{Mod}\big) \simeq H^*_{dR}(\mu\inv(0)^{ss}/G)$.
\end{prop}

Finally, we remark that the identifications of Propositions \ref{HH of I} and \ref{HH of DQ} both come with shifts: $HH_k \simeq H_{dR}^{2n-k}$, where 
$2n = \on{dim}(T^*(X/G)) = \on{dim}(\mu\inv(0)^{ss}/G)$; the latter identification of dimensions comes since the moment map $\mu$ is assumed to be flat, hence $T^*(X/G)$ is equidimensional.

\bibliographystyle{alpha}

\end{document}